\documentclass{agtart_a}
\pdfoutput=1

\usepackage{pinlabel}
\usepackage[all]{xy}

\title{Knot Floer homology in cyclic branched covers}

\author{J Elisenda Grigsby}
\givenname{J Elisenda}
\surname{Grigsby}
\address{Evans Hall\\University of California, Berkeley\\\newline
Berkeley, CA 94720\\USA}
\email{egrigsby@math.berkeley.edu}
\urladdr{}

\volumenumber{6}
\issuenumber{}
\publicationyear{2006}
\papernumber{49}
\startpage{1355}
\endpage{1398}

\doi{}
\MR{}
\Zbl{}


\received{9 September 2005}
\revised{}
\accepted{10 June 2006}
\published{25 September 2006}
\publishedonline{25 September 2006}
\proposed{}
\seconded{}
\corresponding{}
\editor{}
\version{}

\arxivreference{math.GT/0507498}

\let\xysavmatrix\xymatrix
\def\xymatrix{\disablesubscriptcorrection\xysavmatrix}
\AtBeginDocument{\let\bar\wbar\let\tilde\wtilde\let\hat\what
}

\makeop{Spin}
\makeop{Eul}
\makeop{sgn}
\makeop{Sym}
\makeop{hd}
\makeop{hb}

\makeatletter
\def\cnewtheorem#1[#2]#3{\newtheorem{#1}{#3}[section]
\expandafter\let\csname c@#1\endcsname\c@theorem}
\makeatother

\newtheorem{theorem}{Theorem}[section]
\cnewtheorem{lemma}[theorem]{Lemma}
\cnewtheorem{proposition}[theorem]{Proposition}
\theoremstyle{definition}
\cnewtheorem{definition}[theorem]{Definition}

\newcommand{\Ozsvath}{{Ozsv{\'a}th  }}
\newcommand{\Szabo}{{Szab{\'o} }}

\begin{document}

\begin{htmlabstract}
In this paper, we introduce a sequence of invariants of a knot K in
S<sup>3</sup>: the knot Floer homology groups
HFK-hat(&Sigma;<sup>m</sup>(K);K<sup>~</sup>,i) of the preimage of
K in the m&ndash;fold cyclic branched cover over K.  We exhibit
HFK-hat(&Sigma;<sup>m</sup>(K);K<sup>~</sup>,i) as the
categorification of a well-defined multiple of the Turaev torsion of
&Sigma;<sup>m</sup>(K) - K<sup>~</sup> in the case where &Sigma;<sup>m</sup>(K)
is a rational homology sphere.  In addition, when K is a two-bridge
knot, we prove that
HFK-hat(&Sigma;<sup>2</sup>(K);K<sup>~</sup>,s<sub>0</sub>) &asymp;
HFK-hat(S<sup>3</sup>;K) for s<sub>0</sub> the spin Spin<sup>c</sup>
structure on &Sigma;<sup>2</sup>(K).  We conclude with a calculation involving
two knots with identical HFK-hat(S<sup>3</sup>;K,i) for which
HFK-hat(&Sigma;<sup>2</sup>(K); K<sup>~</sup>,i) differ as
<b>Z</b><sub>2</sub>&ndash;graded groups.
\end{htmlabstract}

\begin{abstract}
In this paper, we introduce a sequence of invariants of a knot $K$ in
$S^3$: the knot Floer homology groups
\mbox{$\widehat{HFK}(\Sigma^m(K);\widetilde{K},i)$} of the preimage of
$K$ in the $m$--fold cyclic branched cover over $K$.  We exhibit
\mbox{$\widehat{HFK}(\Sigma^m(K);\widetilde{K},i)$} as the
categorification of a well-defined multiple of the Turaev torsion of
\mbox{$\Sigma^m(K) - \widetilde{K}$} in the case where $\Sigma^m(K)$
is a rational homology sphere.  In addition, when $K$ is a two-bridge
knot, we prove that
\mbox{$\widehat{HFK}(\Sigma^2(K);\widetilde{K},\mathfrak{s}_0) \cong
\widehat{HFK}(S^3;K)$} for $\mathfrak{s}_0$ the spin Spin$^c$
structure on $\Sigma^2(K)$.  We conclude with a calculation involving
two knots with identical \mbox{$\widehat{HFK}(S^3;K,i)$} for which
\mbox{$\widehat{HFK}(\Sigma^2(K); \widetilde{K},i)$} differ as
$\mathbb{Z}_2$--graded groups.
\end{abstract}

\begin{asciiabstract}
In this paper, we introduce a sequence of invariants of a knot K in
S^3: the knot Floer homology groups hat{HFK}(Sigma^m(K);tilde{K},i) of
the preimage of K in the m-fold cyclic branched cover over K.  We
exhibit hat{HFK}(Sigma^m(K);tilde{K},i)$ as the categorification of a
well-defined multiple of the Turaev torsion of Sigma^m(K) - tilde{K}
in the case where Sigma^m(K) is a rational homology sphere.  In
addition, when K is a two-bridge knot, we prove that
hat{HFK}(Sigma^2(K);tilde{K},s_0) is isomorphic to hat{HFK}(S^3;K) for
s_0 the spin Spin^c structure on Sigma^2(K).  We conclude with a
calculation involving two knots with identical hat{HFK}(S^3;K,i)$ for
which hat{HFK}(Sigma^2(K);tilde{K},i) differ as Z_2-graded groups.
\end{asciiabstract}

\maketitle  

\section{Introduction}

Let $Y$ be a closed, connected, oriented 3--manifold and $\mathfrak{s}$ a Spin$^c$ structure on $Y$.  In \cite{MR2113019}, \Ozsvath and \Szabo assign to the pair $(Y; \mathfrak{s})$ a graded abelian group, denoted $\widehat{HF}(Y; \mathfrak{s})$. 

The additional data of an oriented, nullhomologous link $L$ in $Y$ induces a filtration on the chain complex used to compute $\widehat{HF}(Y; \mathfrak{s})$ for each Spin$^c$ structure $\mathfrak{s}$ \cite{MR2065507}, \cite{GT0306378} .  The filtered chain homotopy type of this complex is an invariant of the oriented link $L$ in $Y$.  One can, in particular, calculate the associated graded object of this filtration, yielding a sequence of graded abelian groups $\widehat{HFK}(Y; L, \mathfrak{s}, i)$, called the knot Floer homology groups of $L$ in $Y$.

Now consider $\Sigma^m(K)$, the $m$--fold cyclic branched cover of $S^3$ branched along $K$.  Let $p\co \Sigma^m(K) \rightarrow S^3$ denote the associated projection map and $\widetilde{K} = p^{-1}(K)$ denote the preimage of $K$ in $\Sigma^m(K)$.  Consideration of $\widetilde{K}$ in each cyclic branched cover, $\Sigma^m(K)$, yields a sequence of invariants of the original knot $K$ in $S^3$.  Namely, for each $m \in \mathbb{Z}^+$ we have: 

\begin{definition} \label{definition:inv}
  $\widehat{HFK}(\Sigma^m(K);\widetilde{K}) = \bigoplus_{\mathfrak{s}, i} \widehat{HFK}(\Sigma^m(K); \widetilde{K}, \mathfrak{s}, i)$, the knot Floer homology groups of $\widetilde{K} \subset \Sigma^m(K)$.
\end{definition}

Our aim here is to study this sequence of invariants, focusing on the
case where $K$ is a two-bridge knot and $m=2$.  Then $\Sigma^2(K)$ is
a lens space (Chapter 12 in \cite{MR1959408}) with $H_1(\Sigma^2(K))
= \mathbb{Z}_n$, $n$ an odd integer.  Our main result, stated more
precisely in \fullref{section:twobridge}, is:

\medskip
{\bf\fullref{theorem:central}}\qua
\sl
For $K$ a two-bridge knot in $S^3$, there exists a Spin$^c$ structure, $\mathfrak{s}_0$, on $\Sigma^2(K)$ for which $$\widehat{HFK}(\Sigma^2(K); \widetilde{K};\mathfrak{s}_0) \cong \widehat{HFK}(S^3;K).$$
\rm

This result falsely suggests that the groups $\widehat{HFK}(\Sigma^m(K); \widetilde{K})$ contain no more information than the groups $\widehat{HFK}(S^3; K)$.  In fact, there are knot pairs $K_1$, $K_2$ in $S^3$ for which $$\widehat{HFK}(S^3; K_1) \cong \widehat{HFK}(S^3; K_2)$$ but for which $$\widehat{HFK}(\Sigma^2(K); \widetilde{K_1}) \not\cong \widehat{HFK}(\Sigma^2(K); \widetilde{K_2})$$ as $\mathbb{Z}_2$--graded groups.  Such a pair (the two-bridge knots $K(15,4)$ and $K(15,7)$) is discussed in detail in \fullref{section:examples}.

We also show that for $K$ a nullhomologous knot in a rational homology sphere $Y$, $\widehat{HFK}(Y;K)$ is a categorification of a multiple of the Turaev torsion of $Y-K$.  The connection, established by Kirk and Livingston in \cite{MR1670420}, between the Casson--Gordon invariant of $K$ and various torsions of $\Sigma^2(K)$ bears further examination, particularly since it may yield new obstructions to $K$ being slice.

The paper is laid out as follows:

In \fullref{section:Handlebody} we recall the relevant definitions and theorems in Heegaard Floer homology as well as describe and develop notation for certain natural handlebody decompositions and Heegaard diagrams associated to $K \subset S^3$ and $\widetilde{K} \subset \Sigma^m(K)$.

In \fullref{section:AP} we discuss torsions of chain complexes and prove that $\widehat{HFK}(Y;K)$ is the categorification of a multiple of the Turaev torsion of $Y-K$ in the case where $Y$ is a rational homology sphere.

In \fullref{section:twobridge} we study the invariant $\widehat{HFK}(\Sigma^2(K); \widetilde{K})$ for the case where $K$ is a two-bridge knot in $S^3$.  We also compute $\widehat{HFK}(\Sigma^2(K); \widetilde{K})$ in a few Spin$^c$ structures for the two-bridge knots $K(15,7)$ and $K(15,4)$, whose double branched covers are the lens spaces $-L(15,7)$ and $-L(15,4)$, respectively.

\subsubsection*{Acknowledgments} I am grateful to Peter \Ozsvath for recommending this problem and for numerous indispensable conversations, to Rob Kirby for his guidance and support and for a careful reading of the manuscript, and to the referee for many useful comments.

\section{Background and conventions} \label{section:Handlebody}

We begin by reminding the reader of the Floer homology setup for $\widehat{HFK}(Y)$.  For details, see \Ozsvath and \Szabo \cite{MR2065507,MR2113020,MR2113019}.  For the knot Floer homology refinements, see also Rasmussen \cite{GT0306378}.

\subsection{Heegaard Floer homology background}
Let $K$ be a nullhomologous knot in a closed, oriented, connected 3--manifold $Y$.  Although \Ozsvath and Szab{\'o}'s theory assigns homology groups more generally to nullhomologous links in $Y$, we will focus on knots in this paper.

In \cite{MR2065507}, \Ozsvath and \Szabo present the data of a knot $K$ in $Y$ by means of a {\it doubly-pointed Heegaard diagram} compatible with $K$.  More specifically, they construct a handlebody decomposition of $Y$ arising from a generic self-indexing Morse function $$f\co Y \rightarrow \mathbb{R}$$ with a single index $0$ and $3$ critical point  and $g$ index $1$ and $2$ critical points.  This decomposition yields a Heegaard diagram for $Y$.  The data of two points on the Heegaard surface, $S$, specifies the knot, $K$.

\begin{definition} \label{definition:dopointHD}
  A doubly-pointed Heegaard diagram for a pair $(Y,K)$ is a tuple $(S, \vec{\alpha}, \vec{\beta}, w, z)$ where 
\begin{itemize}
  \item $\vec{\alpha}$ is the $g$--tuple of co-attaching circles for the $g$ $1$--handles
  \item $\vec{\beta}$ is the $g$--tuple of attaching circles for the $g$ $2$--handles
  \item $w, z \in S - \vec{\alpha} - \vec{\beta}$ 
  \item $K$ is the isotopy class of $-\gamma_w \cup \gamma_z$, where $\gamma_w$ and $\gamma_z$ are gradient flow lines from the index $3$ to index $0$ critical points associated to any generic metric on $Y$, intersecting $S$ at $z$ and $w$, respectively.
\end{itemize}
\end{definition}

We gather the standard definitions and notation here for the reader's convenience:

\begin{itemize}
\item $\Sym^g(S) = S^{\times g}/\Sigma_g$ is the $g$--fold symmetric product of the Heegaard surface, $S$.
\item $\mathbb{T}_\alpha = \alpha_1 \times \ldots \times \alpha_g$ (resp.\ $\mathbb{T}_\beta$) is the half-dimensional torus of co-attaching (resp.\ attaching) circles of the $1$--handles (resp.\ $2$--handles) inside $\Sym^g(S)$.
\item $V_z = \{z\} \times \Sym^{g-1}(S)$ (resp.\ $V_w$) is the codimension$_{\mathbb{C}}$ $1$ subvariety of $\Sym^g(S)$ consisting of $g$--tuples where one point is constrained to lie at $z$ (resp.\ $w$).
\item $n_z = (\, - \,) \cap V_z$ (resp.\ $n_w = (\, - \,) \cap V_w$) is the algebraic intersection number of a class with $V_z$ (resp.\ $V_w$) in $\Sym^g(S)$.\end{itemize}

This data gives rise to a $\mathbb{Z}$--filtered chain complex, $\widehat{CF}$ whose 

\begin{itemize}
  \item
    generators are elements $\begin{bf} x \end{bf} \in \mathbb{T}_\alpha \cap \mathbb{T}_\beta$,
  \item
    differential is given by:  $$\widehat{\partial}\begin{bf}x \end{bf} = \sum_{\begin{bf}y \end{bf} \in \mathbb{T}_\alpha \cap \mathbb{T}_\beta} \sum_{\{\phi \in \pi_2(\begin{bf}x, y \end{bf})|\mu(\phi)=1, n_w(\phi) = 0\}} \#(\widehat{\mathcal{M}}(\phi))\,\, \begin{bf}y \end{bf}$$
  \item
where $\widehat{\mathcal{M}}(\phi)$ is the moduli space of holomorphic maps of the standard unit disk into $\Sym^g(S)$ with boundary conditions as in \cite{MR2113019}, in the homotopy class of $\phi$, modded out by the standard $\mathbb{R}$ action (for the sake of simplicity, count $\# \widehat{\mathcal{M}}(\phi)$ with $\mathbb{Z}_2$ coefficients), 
  \item $\mu(\phi)$ is the expected dimension of the moduli space $\mathcal{M}(\phi)$ (before we mod out by the $\mathbb{R}$ action), given by the Maslov index of $\phi$,
  \item and the relative $\mathbb{Z}$ filtration on generators is given by $$\mathcal{F}(\begin{bf} x \end{bf}) - \mathcal{F}(\begin{bf}y\end{bf}) = n_z(\phi) - n_w(\phi)$$ where $\phi \in \pi_2(\begin{bf}x,y\end{bf})$.
\end{itemize}

The relative $\mathbb{Z}$ filtration is improved to an absolute $\mathbb{Z}$ filtration by requiring that the Euler characteristic of the associated bi-graded complex is the symmetrized Alexander polynomial of $K$.  We will have more to say on this point in \fullref{section:AP}.  See also Section 2.3 of \cite{MR2065507}.

The homology groups of the associated graded object of this $\mathbb{Z}$--filtered complex are \Ozsvath and  Szab{\'o}'s knot Floer homology groups; ie, $\widehat{HFK}(Y;K,j)$ is the homology of the chain complex $\mathcal{F}_j/\mathcal{F}_{j-1}$.

For convenience, we introduce a couple more pieces of notation:

\begin{itemize}
  \item $s({\bf x}, {\bf y})= n_z(\phi) - n_w(\phi)$ is the $\mathbb{Z}$--filtration difference between ${\bf x}$ and ${\bf y}$ in $\widehat{CFK}$ (where $\phi \in \pi_2({\bf x}, {\bf y})$).
  \item $m({\bf x}, {\bf y}) = \mu(\phi) - 2n_w(\phi)$ is the relative homological grading of ${\bf x}$ and ${\bf y}$ in $\widehat{CFK}$ (where, again, $\phi \in \pi_2({\bf x}, {\bf y})$).
\end{itemize}

\subsection{Natural handlebody decompositions} \label{subsection:HBody}
In what follows, whenever we refer to a handlebody decomposition for $S^3 - K$, we will always mean one with

\begin{itemize}
  \item a single $0$--handle $h_0$, 
  \item $g$ $1$--handles $h_{\alpha_1}, \ldots h_{\alpha_g}$, 
  \item ($g-1$) $2$--handles $h_{\beta_1}, \ldots, h_{\beta_{g-1}}$, 
  \item no $3$--handles.
\end{itemize}

We will also specify an oriented meridian, $\mu$, for $K$ (along which the final $2$--handle, $h_{\beta_g}$, will be attached to build $S^3$) such that the attaching circle of $h_{\beta_g}$ goes over one of the $h_{\alpha_i}$ (for definiteness, $h_{\alpha_g}$) geometrically once and over all of the other $1$--handles geometrically $0$ times.  $\mu$ generates $H_1(S^3 - K)$ and specifies an orientation.  

We will use the notation $\hb(S^3-K;\mu)$ to denote such a handlebody decomposition.  

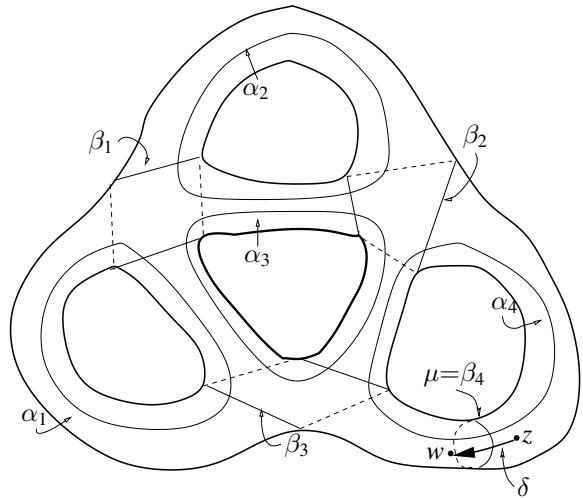
\begin{figure}[ht!]
\begin{center}
\labellist
\pinlabel $\alpha_2$ [t] at 152 251
\pinlabel $\alpha_3$ [t] at 153 148
\pinlabel $\alpha_1$ [b] at 16 36
\pinlabel $\alpha_4$ [b] at 305 105
\pinlabel $\beta_1$ [r] at 67 212
\pinlabel $\beta_2$ [b] at 287 206
\pinlabel $\beta_3$ [l] at 167 26
\pinlabel $\mu{=}\beta_4$ [b] at 274 60
\pinlabel $w$ [r] at 270 21
\pinlabel $z$ [l] at 313 31
\pinlabel $\delta$ [l] at 310 0
\endlabellist
\includegraphics[width=3in]{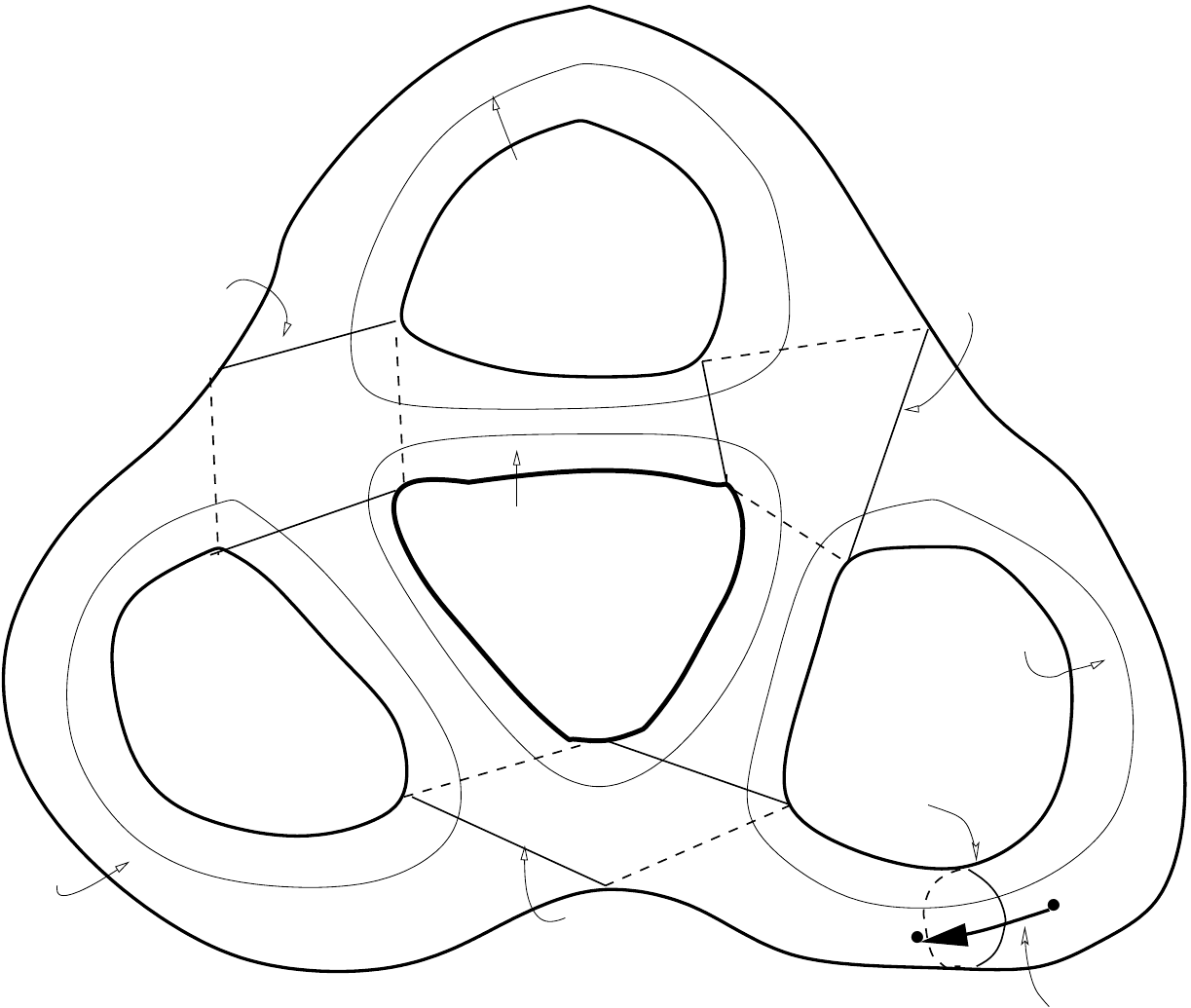}
\end{center}
\caption{Example of a doubly-pointed Heegaard diagram}
\label{fig:dopointHD}
\end{figure}

Similarly, $\hb(S^3;K;\mu)$ will denote the extension of $\hb(S^3 - K; \mu)$ to a handlebody decomposition for $S^3$.  In particular, $$\hb(S^3;K;\mu) = \hb(S^3 - K;\mu) \cup h_{\mu} \cup h_3.$$
Accordingly, we construct a doubly-pointed Heegaard diagram $\hd(S^3;K;\mu)$ for the pair $(S^3,K)$ by choosing an oriented arc $\delta$ on $S$ meeting $\mu$ transversely in a single intersection point and having $0$ geometric intersection with all other co-attaching (attaching) circles for the $1$--handles ($2$--handles).  Our two basepoints $z$ and $w$ are then the initial and final points, respectively, of $\delta$.  See \fullref{fig:dopointHD}.

$\hd(S^3;K;\mu)$ has the properties:

\begin{itemize}
  \item $\vec{\alpha} = \alpha_1 \cup \ldots \cup \alpha_g$ are the coattaching circles associated to the $1$--handles, $h_{\alpha_1}, \ldots h_{\alpha_g}$,
  \item $\vec{\beta} = \beta_1 \cup \ldots \cup \beta_{g-1} \cup (\beta_g = \mu)$ are the attaching circles associated to the $2$--handles $h_{\beta_1}, \ldots, h_{\beta_{g-1}}, h_{\beta_g = \mu}$,
  \item The orientation convention for $K$ given in \fullref{definition:dopointHD} has the property that if $\lambda \subset S$ is a longitude which agrees with the orientation on $K$, then $\lambda \cap \mu = \delta \cap \mu$.
\end{itemize}

A natural handlebody decompostion and doubly-pointed Heegaard diagram can be constructed for the $m$--fold cyclic branched cover $\Sigma^m(K)$ as follows:  
\begin{enumerate}
\item 
  Begin with the natural $\mathbb{Z}_m$--equivariant handlebody decomposition of $\Sigma^m(K)-\widetilde{K}$ which associates to each handle, $h$, in the handlebody decomposition of $S^3-K$, $m$ handles in the handlebody decomposition of $\Sigma^m(K) - \widetilde{K}$ consisting of a chosen lift $\widetilde{h}$ of $h$ and $m-1$ translates $\tau_m(\widetilde{h}), \tau_m^2(\widetilde{h}), \ldots, \tau_m^{m-1}(\widetilde{h})$ of $\widetilde{h}$.

We now have a $\mathbb{Z}_m$--equivariant handlebody decomposition for \sloppy\mbox{$\Sigma^m(K) - \widetilde{K}$} but too many $0$--handles (the theory requires a handlebody decomposition for $\Sigma^m(K)$ with a single $0$ and $3$ handle)\fussy.  

To correct this, recall that $\mu$ is the core circle for a single one of the $1$--handles, $h_{\alpha_g}$.  Use $m-1$ of the lifts of $h_{\alpha_g}$: $\tau_m(\widetilde{h}_{\alpha_g}), \ldots, \tau_m^{m-1}(\widetilde{h}_{\alpha_g})$ to cancel the extra $0$--handles $\tau_m(\widetilde{h}_0), \ldots, \tau_m^{m-1}(\widetilde{h}_0)$.  This new handlebody decomposition is still $\mathbb{Z}_m$--equivariant with respect to the projection map (only now the action on $\widetilde{h}_0$ and $\widetilde{h}_{\alpha_g}$ are trivial).  We denote this handlebody decomposition by $\hb(\Sigma^m(K) - \widetilde{K};\widetilde{\mu})$.
\item
  We extend this to a $\mathbb{Z}_m$--equivariant handlebody decomposition of $\Sigma^m(K)$ by adding one more $2$--handle attached along $\widetilde{\mu}$, which will be our choice of meridian for the knot $\widetilde{K}$ in $\Sigma^m(K)$, and a $3$--handle to fill in the rest of the solid torus neighborhood of $\widetilde{K}$.  We denote this handlebody decomposition by $\hb(\Sigma^m(K); \widetilde{K};\widetilde{\mu})$.
\item
  Associated to this $\hb(\Sigma^m(K); \widetilde{K}; \widetilde{\mu})$ is the corresponding doubly-pointed Heegaard diagram $\hd(\Sigma^m(K); \widetilde{K}; \widetilde{\mu})$ with basepoints $w$ and $z$ on either side of $\widetilde{\mu}$.
\end{enumerate}

An example should make everything concrete.

\subsection{Example: $K$ = right handed trefoil}

\begin{figure}[ht!]
\begin{center}
\labellist
\pinlabel $\alpha_2$ [b] at 36 130
\pinlabel $\alpha_1$ [t] at 279 69
\pinlabel $\beta_1$ [b] at 175 202
\pinlabel $\beta_2{=}\mu$ [l] at 138 121
\endlabellist
\includegraphics[width=3in]{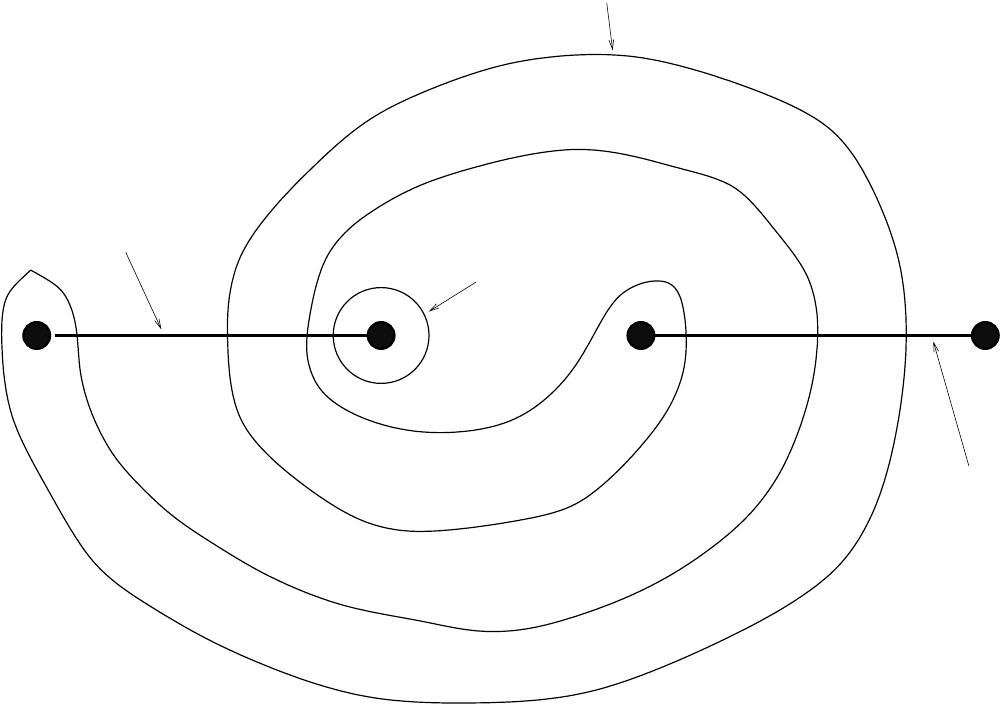}
\end{center}
\caption{Genus $2$ Heegaard diagram for the RH Trefoil}
\label{fig:genus2trefoil}
\end{figure}
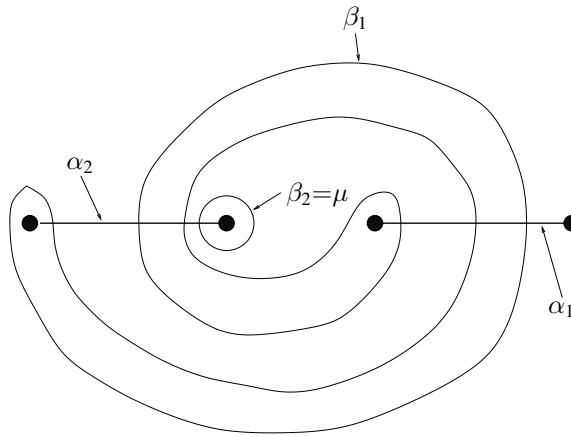

Consider the genus $2$ Heegaard diagram for $S^3$ compatible with $K$ = the right handed trefoil pictured in \fullref{fig:genus2trefoil}.  This is what we have been calling $\hd(S^3;K; \mu = \beta_2)$.  Notice that if we remove the $3$--handle and $h_{\beta_2}$ we get a handlebody decomposition for $S^3 - K$.

A presentation for $\pi_1(S^3 - K)$ is generated by $h_{\alpha_1}, h_{\alpha_2}$ with the single relation given by the attaching map of the $2$--handle $h_{\beta_1}$.  This relation can be read off by traversing $\beta_1$ and keeping track of the intersections with $\alpha_1$ and $\alpha_2$:
$$\pi_1(S^3 - K) = \langle h_{\alpha_1}, h_{\alpha_2}| \partial(h_{\beta_1}) = h_{\alpha_1}h_{\alpha_2}h_{\alpha_1}h_{\alpha_1}^{-1}h_{\alpha_2}^{-1}h_{\alpha_1}^{-1}\rangle$$
The associated Heegaard diagram for $\Sigma^2(K)$ that we have been calling \sloppy $\hd(\Sigma^2(K); \widetilde{K}; \widetilde{\beta}_2)$ has the following properties:

\begin{itemize}
  \item
    $\widetilde{\alpha}_1, \tau_2(\widetilde{\alpha}_1), \widetilde{\alpha}_2$ are the co-attaching circles for the lifts of the $1$--handles (note that we have used $\tau_2(\widetilde{\alpha}_2)$ to cancel the extra $0$--handle)
  \item
    $\widetilde{\beta}_1, \tau_2(\widetilde{\beta}_1), \widetilde{\beta}_2$ are the attaching circles for the lifts of the $2$--handles (note that $\widetilde{\beta}_2$ still intersects $\widetilde{\alpha}_2$ geometrically once and intersects no other $\alpha$ curves)
\end{itemize}

One obtains a handlebody decomposition for $\Sigma^2(K) - \widetilde{K}$ by removing the $3$--handle and the $2$--handle, $\widetilde{h}_{\beta_2}$.

A presentation for $\pi_1(\Sigma^2(K) - \widetilde{K})$ is therefore generated by $\widetilde{h}_{\alpha_1}, \widetilde{h}_{\tau_2(\alpha_1)}, \widetilde{h}_{\alpha_2}$ with the following relations:
\begin{eqnarray*}
  \partial(\widetilde{h}_{\beta_1}) &=& \widetilde{h}_{\alpha_1} \cdot \widetilde{h}_{\alpha_1} \cdot \widetilde{h}_{\alpha_2}^{-1} \cdot \tau_2(\widetilde{h}_{\alpha_1})^{-1} \cdot \widetilde{h}_{\alpha_2}^{-1}\\
  \partial(\tau_2(\widetilde{h}_{\beta_1})) &=& \tau_2(\widetilde{h}_{\alpha_1}) \cdot \widetilde{h}_{\alpha_2} \cdot \tau_2(\widetilde{h}_{\alpha_1}) \cdot \widetilde{h}_{\alpha_1}^{-1}
\end{eqnarray*} 

\section[HFK(Y) and tau(Y-K)]{$\widehat{HFK}(Y)$ and $\check{\tau}(Y-K)$} \label{section:AP}

The aim of this section is to understand $\widehat{HFK}(Y;K)$ as the categorification of a well-defined multiple of a version of the Turaev torsion of $Y-K$ in the case where $Y$ is a rational homology sphere.  $K$, as before, is a nullhomologous knot in $Y$. 

To understand $\widehat{HFK}(Y;K)$ as a categorification, we must
first recall that the chain complex for $\widehat{CFK}(Y;K)$ splits as
a sum of chain complexes, naturally indexed by elements of
Spin$^c(Y_0(K))$, where $Y_0(K)$ denotes the canonical $0$ surgery on
$K$.  The set, Spin$^c(Y_0(K))$, is often referred to as the set of
relative Spin$^c$ structures of the pair $(Y,K)$ and denoted
{\underline{Spin}}$^c(Y,K)$.  In Section 2.3 of \cite{MR2065507} (see
also Section 2.6 of \cite{MR2113019}), \Ozsvath and \Szabo describe,
given a doubly-pointed Heegaard diagram, how to construct a map
$$\mathbb{T}_{\alpha} \cap \mathbb{T}_{\beta} \rightarrow
{\underline{\Spin}}^c(Y,K)$$ and a splitting
$${\underline{\Spin}}^c(Y,K) \cong \Spin^c(Y) \times \mathbb{Z}.$$
Taking the Euler characteristic of each summand in this splitting will
yield a formal polynomial in two variables, one of which indexes the
Spin$^c(Y)$ structure and one of which indexes the $\mathbb{Z}$
factor.

In brief, we obtain such a formal two-variable polynomial as a multiple of the Reidemeister torsion of the maximal abelian cover of $Y-K$.  Specifically, computation of the Reidemeister torsion yields a rational expression in $\mathbb{Q}(H_1(Y))(T)$ which depends on the choice of a lift of a $\mathbb{Z}$--basis for $C_*(Y-K)$ to a $\mathbb{Z}[H_1(Y)][T,T^{-1}]$--basis for $C_*(Y-K)^{[H_1(Y-K)]}$, where $(Y-K)^{[H_1(Y-K)]}$ is the maximal abelian cover of $Y-K$.  A multiple of this rational expression yields a well-defined formal two-variable polynomial in $\mathbb{Z}[\Spin^c(Y)][T,T^{-1}]$ once we use Turaev's correspondence between Spin$^c$ structures on 3--manifolds and lifts of $\mathbb{Z}$--module bases to $\mathbb{Z}[H_1(Y)]$ bases of the maximal abelian cover of $Y$.

Before launching into a formal discussion of these ideas, we state our main result.

\begin{theorem} \label{theorem:categorification}
Let $Y$ be an oriented rational homology sphere and $K$ an oriented, nullhomologous knot in $Y$.

Let $e\co \pi_1(Y-K) \rightarrow \mathbb{Z}$ be given by $e(\gamma) = lk(\gamma, K)$.  Let $r\co \pi_1(Y-K) \rightarrow H_1(Y)$ be the projection onto $H_1(Y)$.

Let $\check{\tau}(Y-K) \in \mathbb{Q}(\Spin^c(Y)(T))$ be the variant of the Reidemeister torsion of the knot complement given in \fullref{definition:Ttorsionvar}.
Then \begin{displaymath}
(\check{\tau}(Y-K)) \cdot (T - 1) = \sum_{\mathfrak{s} \in \mbox{Spin}^c(Y)} p_{\mathfrak{s}}(T) \cdot \mathfrak{s}
\end{displaymath} where
\begin{displaymath}
p_{\mathfrak{s}}(T) = \sum_i \chi(\widehat{HFK}(Y;K,\mathfrak{s},i)) \cdot T^i.
\end{displaymath} 
Recall that $$\chi(\widehat{HFK}(Y; K, \mathfrak{s},i)) = \sum_{\{d \in d_0 + \mathbb{Z}\}} (-1)^{d-d_0} rk(\widehat{HFK}_d(Y;K,\mathfrak{s},i)).$$
\end{theorem}

In the above, $d$ is the absolute homological grading of a generator, defined in \cite{MR2031164}.  We can, however, define $\chi(\widehat{HFK}(Y;K, \mathfrak{s},i))$ without reference to this absolute grading for a rational homology sphere by using the relative $\mathbb{Z}_2$ homological grading on generators induced by comparing the local intersection numbers of $\mathbb{T}_{\alpha}$ and $\mathbb{T}_{\beta}$ at two generators ${\bf x}$ and ${\bf y}$ (see \fullref{theorem:localint}).  

We then lift this relative $\mathbb{Z}_2$ grading to an absolute $\mathbb{Z}_2$ grading by making the choice which insures $$\sum_j (-1)^j rk(\widehat{HF}_j(Y)) = |H_1(Y; \mathbb{Z})|$$ for $j \in \mathbb{Z}_2$.

One should think of $\check{\tau}(Y-K)$ as a rational function in the formal variables $\{\mathfrak{s}|\mathfrak{s} \in$ Spin$^c(Y)\}$ and $T$.  In the numerator of a particular form of this rational expression, the $\mathfrak{s}$ term records the Spin$^c$ structure and the exponent on the $T$ variable records the filtration level of a generator in $\widehat{HFK}(Y)$.

\subsection{Background on torsions of chain complexes}

We start by recalling some definitions.  The classical references for this material are the papers of Milnor \cite{MR0196736}, \cite{MR0242163}.  The particular version of interest to us is developed in Turaev \cite{MR1484699}.

First, recall that the torsion can be defined for a finite acyclic chain complex 
\begin{displaymath}
\xymatrix{0 \ar[r] & C_m \ar[r]^{\partial_m} & \ldots \ar[r]^{\partial_2} & C_1 \ar[r]^{\partial_1} & C_0 \ar[r]^{\partial_0} & 0}\end{displaymath} of vector spaces over a field $\mathbb{F}$ with fixed $\mathbb{F}$--bases $\{c_q\}$ for each $C_q$.  One then chooses, for each $q$, a collection $\{b_q\}$ of elements whose images form a basis for $im(\partial_{q})$ in $C_{q-1}$.

Given two bases $\{b\}$ and $\{c\}$ for a given vector space, let $[b/c]$ denote the determinant of the change of basis matrix (ie, the nonsingular matrix $A = (a_{ij})$ where $b_i = \sum_{j=1}^n a_{ij}c_j$).

Then 

\begin{definition}
The torsion of $C_*$ is defined as $$\tau(C_*) = \prod_{q=0}^m[\{\partial_{q+1}(b_{q+1}), b_q\}/\{c_q\}]^{(-1)^{q+1}}.$$
\end{definition}

It is well-known that $\tau(C_*)$ depends only upon the original choices of bases $\{c_q\}$ for $C_q$.

Classically, we have been interested in torsions of chain complexes arising as covers.  For example, the Alexander polynomial of a knot is (a multiple of) the torsion of the infinite cyclic cover of the knot complement.

In general, we start with a finite chain complex of $\mathbb{Z}$--modules and construct a cover of $X$ via a surjective homomorphism  $\pi_1(X) \rightarrow G$.  Call such a cover $X^{[G]}$.  $C_*(X^{[G]})$ is a free $\mathbb{Z}[G]$--module with a basis given by a choice of lift of the $\mathbb{Z}$--module basis downstairs.  If $G$ is an abelian group, we can construct the field of fractions $\mathbb{Q}(G)$ of $\mathbb{Z}[G]$ by inverting all non-zerodivisors.  The free $\mathbb{Z}[G]$--module basis we chose for $C_*(X^{[G]})$ then becomes a free $\mathbb{Q}(G)$--basis for $C_*(X^{\mathbb{Q}(G)}) := \mathbb{Q}(G) \otimes_{\mathbb{Z}[G]} C_*(X^{[G]})$.  If $C_*(X^{[G]})$ is an acyclic complex, we can compute its torsion.

\subsection{Torsion of $(Y-K)_{\mathbb{Q}(G \times \mathbb{Z})}$}
Let $K$ be a nullhomologous, oriented knot in an oriented rational homology sphere, $Y$, and let $\pi = \pi_1(Y-K)$.  We are interested in the torsion of the chain complex arising from the surjective Hurewicz homomorphism $$\pi \rightarrow H_1(Y-K).$$ 
Notice that a choice of oriented meridian $\mu$ for the knot yields a splitting $$\varphi_{\mu} = r \times e\co \xymatrix{H_1(Y-K) \ar[r]^{\cong} & H_1(Y) \times \mathbb{Z}} \cong G \times \mathbb{Z}$$ specified by $$\varphi_{\mu}(h) = (h-e(h)\mu,e(h)).\footnote{Note that I am abusing notation here (and will throughout), denoting an element of $H_1$ or $\pi_1$ by a loop representing it, and vice versa.}$$
Here, $$e\co \pi \rightarrow \mathbb{Z}$$ is given by linking number with $K$: $$e(h) = lk(h,K)$$ and $$r\co \pi \rightarrow H_1(Y) \cong G$$ is the composition $\pi_1(Y-K) \rightarrow \pi_1(Y)$ with $\pi_1(Y) \rightarrow H_1(Y)$.
 
Let $\epsilon: \mathbb{Z}[\pi] \rightarrow \mathbb{Z}[\mathbb{Z}]$ denote the $\mathbb{Z}$--ring extension of $e$ and let $\rho\co \mathbb{Z}[\pi] \rightarrow \mathbb{Z}[G]$ denote the $\mathbb{Z}$--ring extension of $r$.

Then let $C_*(\widetilde{Y-K})$ denote the $\mathbb{Z}$--linear chain complex of the universal cover, $\widetilde{Y-K}$, of $Y-K$ and form the free $\mathbb{Z}[G \times \mathbb{Z}]$--module $\mathbb{Z}[G \times \mathbb{Z}] \otimes_{\rho \otimes \epsilon} C_*(\widetilde{Y-K})$ and denote it by $$C_*((Y-K)^{[G \times \mathbb{Z}]}).$$
At this point, we can form $\mathbb{Q}(G \times \mathbb{Z})$, the ring of quotients of $\mathbb{Z}[G \times \mathbb{Z}]$ and construct $$C_{\mathbb{Q}(G \times \mathbb{Z})} := \mathbb{Q}(G \times \mathbb{Z}) \otimes_{\mathbb{Z}[G \times \mathbb{Z}]} C_*((Y-K)^{[G \times \mathbb{Z}]}).$$
Now, to compute $\tau(C_{\mathbb{Q}(G \times \mathbb{Z})})$, we fix a handlebody decomposition for $Y-K$ with 

\begin{itemize} \label{list:hbody}
  \item a single $0$--handle $h_0$,
  \item $g$ $1$--handles $h_{\alpha_1}, \ldots, h_{\alpha_g}$, where, again, we are choosing $h_{\alpha_g}$ to be $\mu$, a meridian, (hence, $e(h_{\alpha_g}) = 1$ and $r(h_{\alpha_g}) = 0$),
  \item $(g-1)$ $2$--handles $h_{\beta_1}, \ldots, h_{\beta_{g-1}}$,
  \item no $3$--handles.
\end{itemize}

which yields a $\mathbb{Z}$--module basis for $C_*(Y-K)$.  

The boundary maps $$\partial_*\co C_*((Y-K)^{[G \times \mathbb{Z}]}) \rightarrow C_{*-1}((Y-K)^{[G \times \mathbb{Z}]})$$ are most easily expressed using Fox calculus (see \cite{MR0140099}).  Specifically $\partial_2$ is the $(g-1) \times g$ matrix
\begin{displaymath}
  (\partial_2)_{ij} = (\rho \otimes \epsilon)\left(\frac{\partial h_{\beta_i}}{\partial h_{\alpha_j}}\right) \in \mathbb{Z}[G \times \mathbb{Z}].
\end{displaymath} 
and $\partial_1$ is the $g \times 1$ matrix
\begin{displaymath}
  (\partial_1)_i = (\rho \otimes \epsilon)(h_{\alpha_i} -1) \in \mathbb{Z}[G \times \mathbb{Z}].
\end{displaymath}
After verifying that the chain complex $C_{\mathbb{Q}(G \times \mathbb{Z})}$ is acyclic (addressed by \fullref{lemma:acyclic}, whose proof we give in \fullref{subsection:acyclic}), we can pick lifts $\widetilde{h}_0, \widetilde{h}_{\alpha_i}, \widetilde{h}_{\beta_j}$ of the the $\mathbb{Z}$--module bases downstairs and compute the torsion by comparing that lift with

\begin{itemize}
  \item $b_2 = \{\widetilde{h}_{\beta_1}, \ldots, \widetilde{h}_{\beta_{g-1}}\}$
  \item $\partial_2(b_2) = \{\partial(\widetilde{h}_{\beta_1}), \ldots, \partial(\widetilde{h}_{\beta_{g-1}})\}$
  \item $b_1 = \{\widetilde{h}_{\alpha_g}\}$
  \item $\partial(b_1) = \{\partial(\widetilde{h}_{\alpha_g})\}$
  \item $b_0 = \{\}$
\end{itemize}

Then,

\begin{itemize}
  \item $[b_2/c_2] = 1$
  \item $[\partial_2(b_2)b_1/c_1] = \rho \otimes \epsilon([\partial_2^g])$, where $\partial_2^g$ is the $g \times g$ matrix obtained by inserting the tuple representing the $g$th basis element of the chosen basis $\{h_{\alpha_1}, \ldots h_{\alpha_g}\}$ into the $g$th row of $\partial_2$.
  \item $[b_1/c_0] = \rho \otimes \epsilon([\partial_1^g])$, where $\partial_1^g$ is the $1 \times 1$ matrix obtained by deleting the first $g-1$ rows from $\partial_1$.  In other words, $[b_1/c_0] = \rho \otimes \epsilon(h_{\alpha_g} - 1) = T^{1} - T^{0}.$
\end{itemize}

We obtain the torsion of $C_{\mathbb{Q}(G \times \mathbb{Z})}$ by computing the rational expression $$(\rho \otimes \epsilon)\left(\frac{[\partial_2^g]}{[\partial_1^g]}\right) = \frac{\rho \otimes \epsilon([\partial_2^g])}{T^{1}-T^{0}}$$ with respect to the initial choices $\widetilde{h}_0, \widetilde{h}_{\alpha_i}, \widetilde{h}_{\beta_j}$ of lifts.

\subsection{Euler chains and bases for chain complexes}

All we lack in the above is a nice way of specifying a $\mathbb{Q}(G \times \mathbb{Z})$--basis for the chain complex $C_{\mathbb{Q}(G \times \mathbb{Z})}$, ie, a lift of a particular $\mathbb{Z}$--module basis for $C_*(Y-K)$.  Without such a lift, $\tau$ has an indeterminacy coming from this choice.

It turns out that Turaev gives us exactly the tools we need to specify such a lift.  In brief, he explains how to associate to an \begin{it} Euler chain \end{it} (defined in Section 2.5 in \cite{MR1484699}) on $Y-K$ a $\mathbb{Z}[H_1(Y-K)]$--module basis for the maximal abelian cover, $(Y-K)^{[H_1(Y-K)]}$.  

Furthermore, given a doubly-pointed Heegaard diagram for $Y$ compatible with a knot $K$, \Ozsvath and \Szabo construct a map $$\mathbb{T}_{\alpha} \cap \mathbb{T}_{\beta} \rightarrow \Eul(Y_0(K)) \rightarrow \Eul(Y-K).$$
There is a natural splitting of $\Eul(Y-K) \cong \Eul(Y) \times \mathbb{Z}$ which, via Turaev's identification of Euler chains and Spin$^c$ structures, allows us to realize the torsion of the chain complex of the maximal abelian cover as a formal element of $\mathbb{Q}(\Spin^c(Y))(T)$.

So, choosing a lift of a $\mathbb{Z}$--module basis for $C_*(Y-K)$ to a \sloppy \mbox{$\mathbb{Z}[H_1(Y-K)]$}-module basis for \mbox{$C_*(Y-K)^{[H_1(Y-K)]}$} \fussy is just a matter of specifying an element of $\Eul(Y-K)$.

We specify such an element by using \Ozsvath and  Szab{\'o}'s identification \cite{MR2065507} $$\mathbb{T}_\alpha \cap \mathbb{T}_\beta \rightarrow \underline{\Spin}^c(Y,K) \leftrightarrow \Eul(Y_0(K)).$$
Their map comes complete with a natural map $$\Eul(Y_0(K)) \rightarrow \Eul(Y-K)$$ induced by forgetting the final two and three handle (and corresponding arcs in the spider-like Euler chain) along with a splitting $$p_1 \times p_2\co \Eul(Y-K) \rightarrow \Eul(Y) \times \mathbb{Z}.$$
Here, the map $p_1\co \Eul(Y-K) \rightarrow \Eul(Y)$ is obtained via the unique extension of an Euler chain for $Y-K$ to one for $Y$.\footnote{Recall that in going from $Y-K$ to $Y$ we add a $2$--handle along a meridian $\mu$ for the knot and a $3$--handle, and $\mu$ has a unique intersection point with a single $\alpha$ curve.  The arc connecting $h_0$ to $h_3$ is uniquely specified by the basepoint, $w$.} 

The map $p_2\co \Eul(Y-K) \rightarrow \mathbb{Z}$ is defined as follows.  Let $\xi \in \Eul(Y-K)$ be the restriction of $\xi' \in \Eul(Y_0(K))$ and $\mathfrak{s}_{\xi'}$ be the element of Spin$^c(Y_0(K))$ associated to $\xi'$ via Turaev's identification.  Then $$p_2(\xi) = \frac{1}{2}\langle c_1(\mathfrak{s}_{\xi'}),[\hat{F}]\rangle,$$ where $[\hat{F}]$ is the homology class of a capped-off Seifert surface for $K$ in $Y_0(K)$.

Notice that this splitting is defined so that it respects the natural $\mu$--induced splitting $$\varphi_\mu\co \xymatrix{H_1(Y-K) \ar[r]^{\cong} & (H_1(Y) \times \mathbb{Z})}$$ specified by $\varphi_\mu(h) = (h -e(h)\mu, e(h))$, where $e\co H_1(Y-K) \rightarrow \mathbb{Z}$ is defined by $$e(h) = lk(h,K) = \langle (h'),[\hat{F}]\rangle,$$ where $h' \in H_1(Y_0(K))$ is the induced image under the inclusion map $$H_1(Y-K) \rightarrow H_1(Y_0(K))$$ (again, we are assuming that $K$ is an oriented knot).

We are finally ready to define the variant of Reidemeister torsion for which $\widehat{HFK}(Y,K)$ is a categorification in the case when $Y$ is a rational homology sphere.

\begin{definition} \label{definition:Ttorsionvar}
  Let $Y$ be a rational homology sphere, $K$ a nullhomologous, oriented knot in $Y$, $\mu$ a choice of meridian for $K$.

Let $\hd(Y;K,\mu)$ be a doubly-pointed Heegaard diagram compatible with $K$.  

Then consider the summands of the formal determinant of the matrix
$$\begin{array}{lll}
& \begin{array}{rll} \hskip 5pt {\beta}_1 & \ldots & {\beta}_{g-1} \end{array} \\
\begin{array}{r} {\alpha}_1 \\ \vdots \\ {\alpha}_{g-1} \end{array} & \left( \begin{array}{cr} \hskip 50pt &  \\  \hskip 50pt &  \\ \hskip 50pt & \end{array} \right) \end{array}
$$
where the entries of the matrix are formal sums of intersection points between the appropriate $\alpha$ and $\beta$ curves, and each intersection point in the matrix is assigned a $\pm 1$ according to its local intersection number (see \fullref{lemma:maslov}).

By acting on this formal sum by \Ozsvath and  Szab{\'o}'s map $$f\co \mathbb{T}_\alpha \cap \mathbb{T}_\beta \rightarrow \underline{\Spin}^c(Y,K)$$ composed with the splitting $$p_1 \times p_2 \co \underline{\Spin}^c(Y,K) \rightarrow \Spin^c(Y) \times \mathbb{Z}$$ we get a formal element of $\mathbb{Z}[\Spin^c(Y)][T,T^{-1}]$; ie, a formal polynomial, $p_{\mathfrak{s}}(T)$, in the variables $\{\mathfrak{s}|\mathfrak{s} \in \Spin^c(Y)\}$ and $T, T^{-1}$.

We define: $$\check{\tau}(Y-K) = \frac{p_{\mathfrak{s}}(T)}{T-1}.$$
\end{definition}

{\bf Remark}\qua
Note that $\check{\tau}(Y-K)$ is actually the Reidemeister torsion of the maximal abelian cover of $Y-K$, where the usual indeterminacy coming from a choice of basis has been eliminated.  

More precisely, let $e\co \pi_1(Y-K) \rightarrow \mathbb{Z}$ be given by $$e(\gamma) = lk(\gamma,K),$$ $$r\co \pi_1(Y-K) \rightarrow H_1(Y)$$ be the projection onto $H_1(Y)$, and $$\epsilon\co \mathbb{Z}[\pi_1(Y-K)] \rightarrow \mathbb{Z}[T,T^{-1}],$$ $$\rho\co \mathbb{Z}[\pi_1(Y-K)] \rightarrow \mathbb{Z}[H_1(Y)]$$ the natural group ring extensions.

Then $$(\rho \otimes
\epsilon)\left(\frac{[\partial_2^g]}{[\partial_1^g]}\right) =
\frac{\rho \otimes \epsilon([\partial_2^g])}{T-1}$$ yields the
Reidemeister torsion associated to the maximal abelian cover of
\sloppy \mbox{$Y-K$}. \fussy If the basis of the chain complex for the
maximal abelian cover is specified by Ozsv{\'a}th--Szab{\'o}'s map
to Spin$^c$ structures, then we arrive at $\check{\tau}(Y-K)$.

$\check{\tau}(Y-K)$ is the variant of Reidemeister torsion which will (when multiplied by $(T-1)$) be the categorification of the knot Floer homology for a rational homology sphere.

\subsection{$C_{\mathbb{Q}(G \times \mathbb{Z})}$ is an acyclic complex} \label{subsection:acyclic}

We return now to the proof of an important point which we left unresolved in an earlier part of this section.

\begin{lemma} \label{lemma:acyclic}
The chain complex $C_{\mathbb{Q}(G \times \mathbb{Z})}$ is acyclic.
\end{lemma}

\proof[Proof of \fullref{lemma:acyclic}]  We need only show that the determinants of the matrices\break $(\rho \otimes \epsilon)(\partial_2^g)$ and $(\rho \otimes \epsilon)(\partial_1^g)$ are units in $\mathbb{Q}(G)(T)$.  

$(\rho \otimes \epsilon)[\partial_1^g] = T-1$ is clearly a nonzerodivisor in $\mathbb{Z}[G][T,T^{-1}]$ and hence a unit in $\mathbb{Q}(G)(T)$.

To see that $[\partial_2^g]$ must be a unit in $\mathbb{Q}(G)(T)$, consider the ring homomorphism $$\varphi\co \mathbb{Q}(G)(T) \rightarrow \mathbb{Q}$$ which sends $T$ and all $h \in G$ to $1$.  If $\varphi(p) = q$ is a unit in $\mathbb{Q}$, then $p$ must be a unit in $\mathbb{Q}(G)(T)$.  

But $\varphi((\rho \otimes \epsilon)(\partial_2^g))$ is exactly the relation matrix for $H_1(Y;\mathbb{Q})$.  $H_1(Y;\mathbb{Q}) = 0$ then implies that $\varphi((\rho \otimes \epsilon)[\partial_2^g])$ is a unit in $\mathbb{Q}$.  \endproof

\subsection[Proof of \ref{theorem:categorification}]{Proof of \fullref{theorem:categorification}} \label{subsection:determinant}

The relationship between generators of the knot Floer homology and summands of the determinant of the differential $$\partial_2\co (C_{\mathbb{Q}(G \times \mathbb{Z})})_2 \rightarrow (C_{\mathbb{Q}(G \times \mathbb{Z})})_1$$ used to compute the Reidemeister torsion is clear, since the generators of $\widehat{CFK}$ are themselves summands of an analogous formal determinant (see, eg, the discussion in the proof of \fullref{proposition:bijection}).

Under this correspondence, the relative element of $\underline{\Spin}^c(Y,K)$ specified by the difference between a pair of generators is an element of \sloppy \mbox{$H_1(Y_0(K)) \cong H_1(Y-K)$.} \fussy  

Furthermore, by the naturality of the $H_1$ action on all of these sets, the splitting $$H_1(Y-K) \rightarrow H_1(Y) \times \mathbb{Z}$$ matches the splitting $$\underline{\Spin}^c(Y,K) \cong \Spin^c(Y) \times \mathbb{Z}.$$
We need only verify that $s({\bf x}, {\bf y})$, the relative filtration grading of two generators,  ${\bf x}$ and ${\bf y}$, defined by $n_z(\phi) - n_w(\phi)$ (for $\phi \in \pi_2({\bf x},{\bf y})$), matches up with the relative $T$ exponent of the corresponding summands and that $m({\bf x},{\bf y})$, the relative Maslov grading of the generators agrees mod $2$ with the relative sign of the summands in the determinant.

\begin{lemma} \label{lemma:filtration}
  Suppose $K$ is an oriented, nullhomologous knot in a closed, connected, oriented 3--manifold $Y$ and $\hd(Y;K,\mu)$ is a doubly-pointed \sloppy \mbox{Heegaard} \fussy diagram as in \fullref{section:Handlebody}.  Then consider the map $e\co \pi \rightarrow \mathbb{Z}$ given by $e(\gamma) = lk(\gamma, K)$.  Let ${\bf x}$, ${\bf y}$ be Floer homology generators and let $\gamma_{{\bf x}}, \gamma_{{\bf y}} \in \pi$ be their corresponding summands in the Fox determinant, $[\partial_2^g]$.  Then $$e(\gamma_{{\bf x}})- e(\gamma_{{\bf y}}) = s({\bf x},{\bf y}).$$
\end{lemma}

\proof[Proof of \fullref{lemma:filtration}]\qua
First, recall that the filtration difference between two generators ${\bf x}$ and ${\bf y}$ is well-defined whenever ${\bf x}$ and ${\bf y}$ are in the same Spin$^c$ structure.  Then there exists a topological disk $\phi \in \pi_2({\bf x},{\bf y})$ and $s(\begin{bf}x,y \end{bf}) = n_z(\phi) - n_w(\phi)$.

Recall also that in $\hd(Y;K)$ we have a distinguished $\beta$ circle, $\mu$, which is  a meridian for $K$, and next to which we place the two basepoints as shown in \fullref{fig:dopointHD}.  Recall (Section 2.13 of \cite{MR2113019}) that we can represent disks in $\pi_2({\bf x}, {\bf y})$ uniquely as $\mathbb{Z}$--linear combinations of \begin{it} fundamental domains \end{it} $\mathcal{D}_i$, which are the closures of the connected components of $S - \alpha_1 - \ldots - \alpha_g - \beta_1 - \ldots - (\beta_g = \mu)$.

Let $\mathcal{D}_z$ and $\mathcal{D}_w$ be the fundamental domains containing the basepoints $z$ and $w$, respectively.  Given a representation of a disk $\phi$ as a $\mathbb{Z}$--linear combination of fundamental domains, $n_z$ and $n_w$ are then the coefficients on $\mathcal{D}_z$ and $\mathcal{D}_w$, respectively.

Now suppose ${\bf x}$ and ${\bf y}$ are two generators in the same Spin$^c$ structure, $\phi \in \pi_2({\bf x}, {\bf y})$ is a disk connecting them, and $\sum_i \mathcal{D}_i$ is the linear combination of fundamental domains representing $\phi$.  Note that $$\partial(\mathcal{D}_w)= -\mu + \mbox{(other stuff)},$$  and $$\partial(\mathcal{D}_z) = \mu +  \mbox{(other stuff)}.$$  Therefore, 
\begin{eqnarray*}
\partial(\phi) &=& n_z(\partial \mathcal{D}_z) + n_w(\partial \mathcal{D}_w) + \mbox{(other stuff)}\\ 
 &=& (n_z - n_w)\mu +  \mbox{(other stuff)}\\
\end{eqnarray*}
Call this (other stuff) $\gamma$.  Note that $\gamma$ is exactly the image under the Hurewicz map $\pi_1(Y-K) \rightarrow H_1(Y-K)$ of $\gamma_{\bf x}^{-1} \cdot \gamma_{\bf y}$, where $\gamma_{\bf x}$ and $\gamma_{\bf y}$ are the summands corresponding to ${\bf x}$ and ${\bf y}$ in the Fox determinant, $[\partial_2^g]$.  In other words, $$e(\gamma) = e(\gamma_{\bf x}^{-1}\cdot \gamma_{\bf y}) = -(e(\gamma_{\bf x}) - e(\gamma_{\bf y})).$$
Furthermore, $$(n_z - n_w)\mu = -\gamma$$ in $H_1(Y-K)$.  Therefore, 
\begin{eqnarray*}
e(\gamma_{\bf x}) - e(\gamma_{\bf y}) &=& -e(\gamma)\\
 &=& (n_z(\phi) - n_w(\phi))\cdot e(\mu)\\
 &=& n_z(\phi) - n_w(\phi)\\
 &=& s({\bf x}, {\bf y}),\\
\end{eqnarray*}
as desired. \endproof

\begin{lemma} \label{lemma:maslov}
  Let ${\bf x}$ and ${\bf y}$ be two elements of $\widehat{CFK}(Y)$ in the same Spin$^c$ structure.  Let $\gamma_{{\bf x}}, \gamma_{\begin{bf}{y}\end{bf}} \in \mathbb{Z}[\pi]$ be the corresponding summands in the Fox determinant, $[\partial_2^g]$.  Then $$(-1)^{m(\begin{bf} x, y \end{bf})} = \sgn(\gamma_{{\bf x}}) \cdot \sgn(\gamma_{{\bf y}}).$$
\end{lemma}

\proof[Proof of \fullref{lemma:maslov}]\qua
Recall the following standard fact from Lagrangian Intersection Floer theory.

\begin{theorem} [Floer, Robbin--Salamon \cite{MR1241874}] \label{theorem:localint}
Let $L_1$ and $L_2$ be two \sloppy \mbox{Lagrangian} \fussy submanifolds in a symplectic manifold $X$. Given $\begin{bf}x, y \end{bf} \in L_1 \cap L_2$ and $\phi \in \pi_2({\bf x}, {\bf y})$ a pseudoholomorphic disk connecting them, we have $$(-1)^{\mu(\phi)} = deg(\begin{bf}x \end{bf}) \cdot deg(\begin{bf}y \end{bf})$$  Here $deg$ denotes the local intersection number of $L_1$ and $L_2$ at ${\bf x}$ and $\mu(\phi)$ is the Maslov index of $\phi$.
\end{theorem}

In our setting, $\mathbb{T}_{\alpha}$ and $\mathbb{T}_{\beta}$ play the role of Lagrangians in the symplectic manifold $\Sym^g(S)$.  Since the mod $2$ Maslov index difference depends only on the local intersection degrees of the two intersection points ${\bf x}$ and ${\bf y}$ in $\mathbb{T}_{\alpha} \cap \mathbb{T}_{\beta}$, we need only prove that $\sgn(\gamma_{{\bf x}}) \cdot \sgn(\gamma_{\begin{bf} y \end{bf}}) = deg(\begin{bf} x \end{bf}) \cdot deg(\begin{bf}y \end{bf})$, a straightforward calculation in local coordinates on $\Sym^g(S)$. \endproof

\subsection{Relationship to twisted Alexander polynomials}

We mention some closely-related constructions developed in Wada 
\cite{MR1273784}, Kirk and Living\-ston \cite{MR1670420} and Kitano \cite{MR1405595}.

As usual, we assume $Y$ is an oriented rational homology sphere, and $K$ in $Y$ is an oriented, nullhomologous knot.

Then we have an isomorphism $$f\co H_1(Y) \rightarrow \mathbb{Z}_{a_1} \times \ldots \times \mathbb{Z}_{a_k}$$ where the elements $a_i$ are well defined under the added condition that $a_i$ divides $a_{i+1}$ for all $1 \leq i < k$.  Let $c_i$ denote $(p_i \circ f)(c)$, where $$p_i\co (\mathbb{Z}_{a_1} \times \ldots \times \mathbb{Z}_{a_k}) \rightarrow \mathbb{Z}_{a_i}$$ denotes the projection onto the $i$th component.  Then for each $i$ we have a character $\phi_i\co H_1(Y) \rightarrow S^1 \subset \mathbb{C}$ which lands in the cyclotomic field $\mathbb{Q}(\zeta_{a_i})$ for $\zeta_{a_i}$ a primitive $a_i$th root of unity.  By multiplying characters, we obtain the $1$--dimensional tensor product representation: 
\begin{eqnarray*} \label{equation:representation} \phi(c) (z) &=& (\phi_1 \otimes \ldots \otimes \phi_k)(c) (z)\\ &=& \zeta_{a_1}^{c_1} \cdots \zeta_{a_k}^{c_k} (z)
\end{eqnarray*}
Now we can form what is known as the $\phi$--twisted Alexander polynomial of $K$ by using $e\co\pi \rightarrow \mathbb{Z}$ coming from the intersection number with a Seifert surface and its $\mathbb{Z}$--ring extension $\epsilon\co \mathbb{Z}[\pi] \rightarrow \mathbb{Z}[\mathbb{Z}]$.

\begin{definition} [Wada]
  Let $Y$, $K$, $\epsilon$ be as above.  Let $\mathbb{F}$ be a field and $\phi\co \mathbb{Z}[\pi] \rightarrow Gl_n(\mathbb{F})$ a representation.  Then the $\phi$--twisted Alexander polynomial of $K$ in $Y$ is the rational expression $$\Delta_{K, \phi}(T) = \left(\frac{(\phi \otimes \epsilon) [\partial_2^g]}{T - 1}\right)$$ where $\partial_2$ is the Fox matrix associated to the presentation  $$\langle h_{\alpha_1}, \ldots h_{\alpha_g}|h_{\beta_1}, \ldots, h_{\beta_{g-1}}\rangle$$ of $\pi_1(Y-K)$ and $\partial_2^g$ is $\partial_2$ with the $g$th column removed.  Here, the presentation is again assumed to have the property that $h_{\alpha_g}$ is a meridian of the knot, implying $e(h_{\alpha_g}) = 1$ and $f(h_{\alpha_g}) = 0.$
\end{definition}

Kitano proves, in \cite{MR1405595}, that Wada's $\phi$--twisted polynomial is the Reidemeister torsion of the chain complex associated to $\phi \otimes \epsilon$.  Since this chain complex is exactly the chain complex $C_{\mathbb{Q}(G \times \mathbb{Z})}$ (the only difference is an extra map $\mathbb{Q}(G) \rightarrow \mathbb{Q}(\zeta_{a_1}, \ldots, \zeta_{a_k}) \subset \mathbb{C}$, yielding an element of $\mathbb{Q}(\zeta_{a_k})[T,T^{-1}]$),\footnote{Note that, since $a_i|a_k$ for all $i$, $\mathbb{Q}(\zeta_{a_1}, \ldots \zeta_{a_k}) = \mathbb{Q}(\zeta_{a_k}).$} Wada's $\phi$--twisted polynomial is the Reidemeister torsion in a slightly different form.

Kirk and Livingston, in \cite{MR1670420}, define yet another version of a $\phi$--twisted Alexander polynomial, which differs slightly from Wada's definition.  They again begin with homomorphisms $e\co \pi \rightarrow \mathbb{Z}$ and $f\co \pi \rightarrow \mathbb{Z}_{a_1} \times \ldots \times \mathbb{Z}_{a_k}$ and form the $\mathbb{Q}(\zeta_{a_k})[T,T^{-1}]$ chain complex $$C_*(Y-K; \mathbb{Q}(\zeta_{a_k})[T,T^{-1}]_\phi) \co= \mathbb{Q}(\zeta_{a_1}, \ldots \zeta_{a_k}) \otimes_{\phi} C_*(\widetilde{Y-K})$$ where here, the $T$ action is given (once a $\mu$ with $e(\mu) = 1$ is chosen) by $$T^n(g \otimes c) = (g \cdot \phi(\mu^{-n})) \otimes (\mu^n \cdot c)$$ extended linearly.

\begin{definition} [Kirk--Livingston]
The $i$th $\phi$--twisted Alexander polynomial, denoted $\Delta_i$, is the order of the torsion of the $i$th homology of $$C_*(Y-K;\mathbb{Q}(\zeta_{a_k})[T,T^{-1}]_\phi),$$ considered as a $\mathbb{Q}(\zeta_{a_k})[T,T^{-1}]$--module.
\end{definition}

They go on to prove that Wada's invariant, labeled $W$, is related to $\Delta_0$ and $\Delta_1$ by the simple formula $$W = \frac{\Delta_1}{\Delta_0}.$$

\section[HFK for double-branched covers of two-bridge knots]{$\widehat{HFK}$ for double-branched covers of two-bridge knots} \label{section:twobridge}

We now turn to exploring $\widehat{HFK}(\Sigma^m(K);\widetilde{K})$ in the case $m=2$ and $K$ a two-bridge knot.  Our main result is \fullref{theorem:central}.

We start by recalling a few standard facts about two-bridge knots.  A good reference is Chapter 12 of \cite{MR1959408}.

First, there is a one-to-one correspondence between isotopy classes of two-bridge knots and lens spaces arising as their double branched covers.

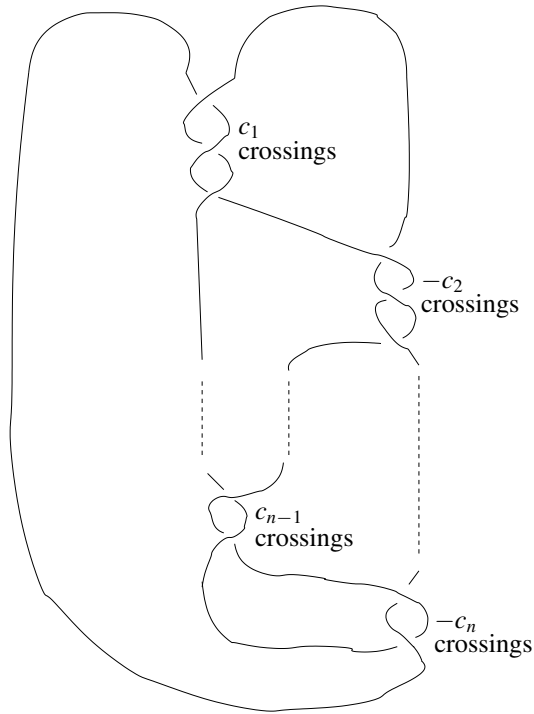
\begin{figure}[ht!]
\begin{center}
\labellist\hair 3pt
\pinlabel $c_1$ [l] at 184 485
\pinlabel crossings [l] <0pt,-10pt> at 184 485
\pinlabel $-c_2$ [l] at 336 358
\pinlabel crossings [l] <0pt,-10pt> at 336 358
\pinlabel $-c_n$ [l] at 347 75
\pinlabel crossings [l] <0pt,-10pt> at 347 75
\pinlabel $c_{n{-}1}$ [l] at 198 163
\pinlabel crossings [l] <0pt,-10pt> at 198 163
\endlabellist
\includegraphics[width=2.2in]{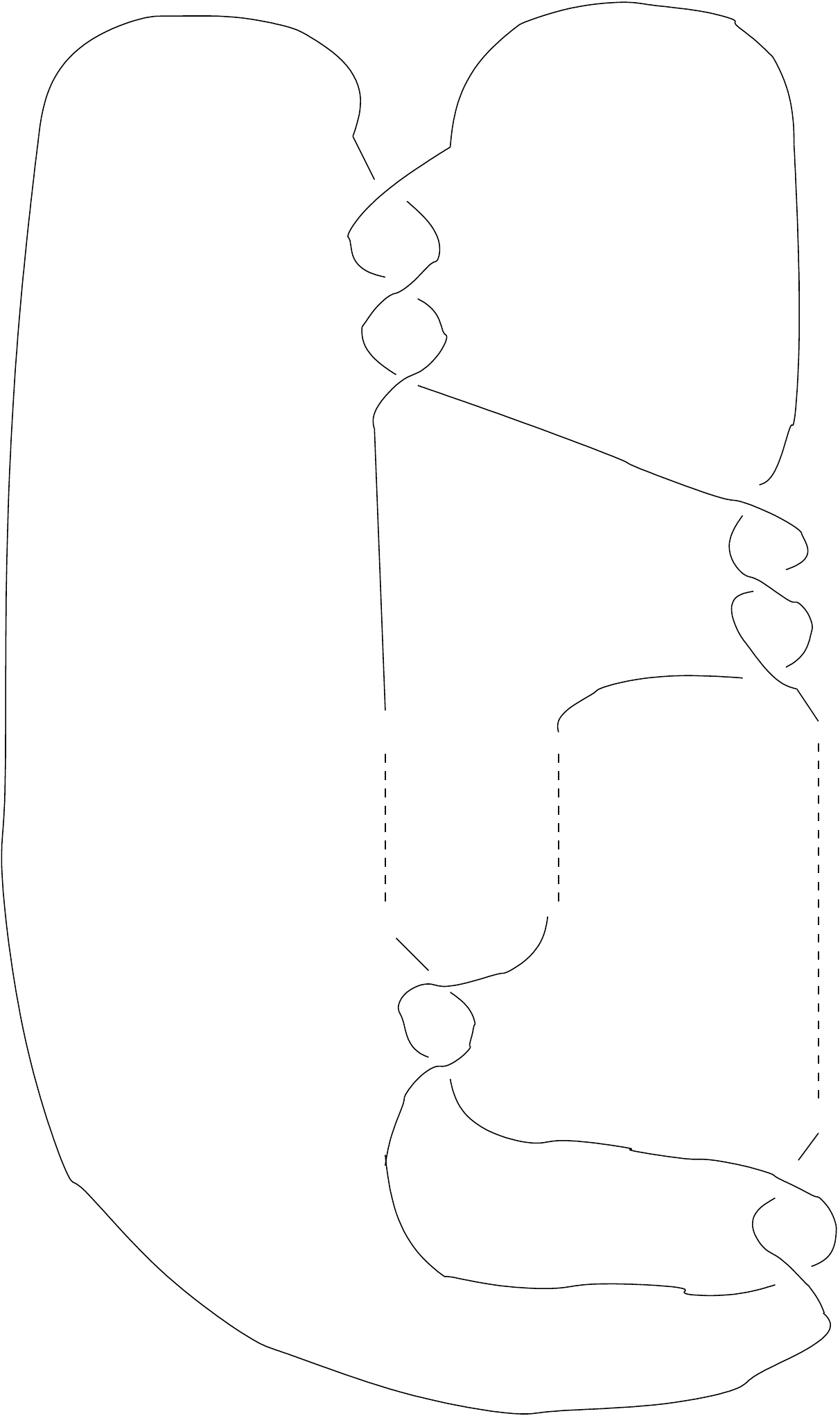}
\end{center}
\caption{A two-bridge knot with twist numbers $(c_1, -c_2, c_3, -c_4, \ldots, c_n)$}
\label{fig:Twobridge}
\end{figure}

\begin{theorem}  {\rm\cite{MR0082104,MR0258014}}\qua A two bridge knot $K$ in $S^3$ with twist numbers $$(c_1, -c_2, c_3, -c_4, \ldots, c_n)$$ (see \fullref{fig:Twobridge}) has double branched covering $-L(p,q)$ where $\frac{p}{q}$ is the continued fraction expansion $$ \frac{p}{q} = c_1 + \frac{1}{\displaystyle c_2 + \frac{1}{\displaystyle \cdots + \frac{1}{c_n}}}.$$
\end{theorem}

We will denote the two-bridge knot whose double branched cover is $-L(p,q)$ by $K(p,q)$.

A particularly useful projection of a two-bridge knot for our purposes is the Schubert normal form.  We construct the Schubert normal form of the knot $K(p,q)$ as a union of 4 segments on $S^2$: 2 straight ``underbridges'' $U_1$ and $U_2$ and two curvy ``overbridges'' $O_1$ and $O_2$ (All of the following is explained very nicely in \cite{MR1928176}).

\begin{figure}[ht!]
\begin{center}
\labellist\hair2pt
\pinlabel $b_p$ [r] at 0 41
\pinlabel $b_0$ [l] at 421 41
\pinlabel $b_{p{-}1}$ [b] at 56 68
\pinlabel $b_5$ [b] <-2pt,0pt> at 120 77
\pinlabel $b_4$ [b] <-2pt,0pt> at 185 80
\pinlabel $b_3$ [b] <-2pt,0pt> at 250 80
\pinlabel $b_2$ [b] <-2pt,0pt> at 314 77
\pinlabel $b_1$ [b] <-2pt,0pt> at 371 67
\pinlabel $b_{p{+}1}$ [t] at 55 13
\pinlabel $b_{2p{-}5}$ [t] at 120 4
\pinlabel $b_{2p{-}4}$ [t] at 185 0
\pinlabel $b_{2p{-}3}$ [t] at 250 0
\pinlabel $b_{2p{-}2}$ [t] at 314 4
\pinlabel $b_{2p{-}1}$ [t] at 371 13
\endlabellist
\includegraphics[width=3in]{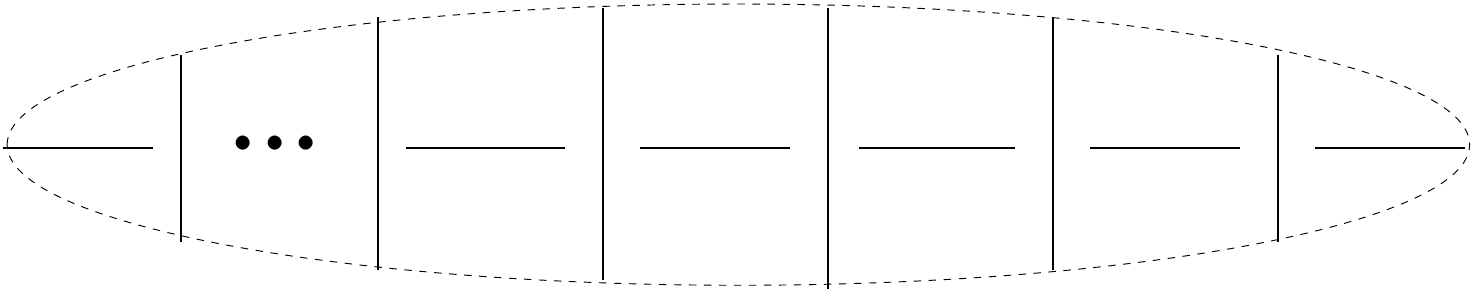}
\end{center}
\caption{A neighborhood of $U_2$}
\label{fig:ubridge2}
\end{figure}

\begin{enumerate}
  \item A neighborhood of $U_1$ looks like \fullref{fig:ubridge1} and a neighborhood of $U_2$ looks like the mirror image of $U_1$ reflected across a central vertical axis as in \fullref{fig:ubridge2}.

\end{enumerate}
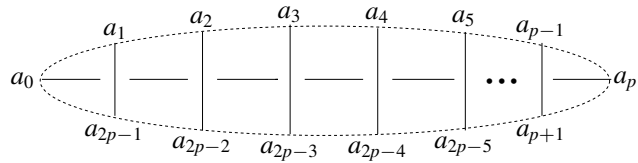
\begin{figure}[ht!]
\begin{center}
\labellist\hair2pt
\pinlabel $a_0$ [r] at 0 41
\pinlabel $a_p$ [l] at 421 41
\pinlabel $a_1$ [b] at 56 68
\pinlabel $a_2$ [b] at 120 77
\pinlabel $a_3$ [b] at 185 80
\pinlabel $a_4$ [b] at 250 80
\pinlabel $a_5$ [b] at 314 77
\pinlabel $a_{p{-}1}$ [b] at 371 67
\pinlabel $a_{2p{-}1}$ [t] at 55 13
\pinlabel $a_{2p{-}2}$ [t] at 120 4
\pinlabel $a_{2p{-}3}$ [t] at 185 0
\pinlabel $a_{2p{-}4}$ [t] at 250 0
\pinlabel $a_{2p{-}5}$ [t] at 314 4
\pinlabel $a_{p{+}1}$ [t] at 371 13
\endlabellist
\includegraphics[width=3in]{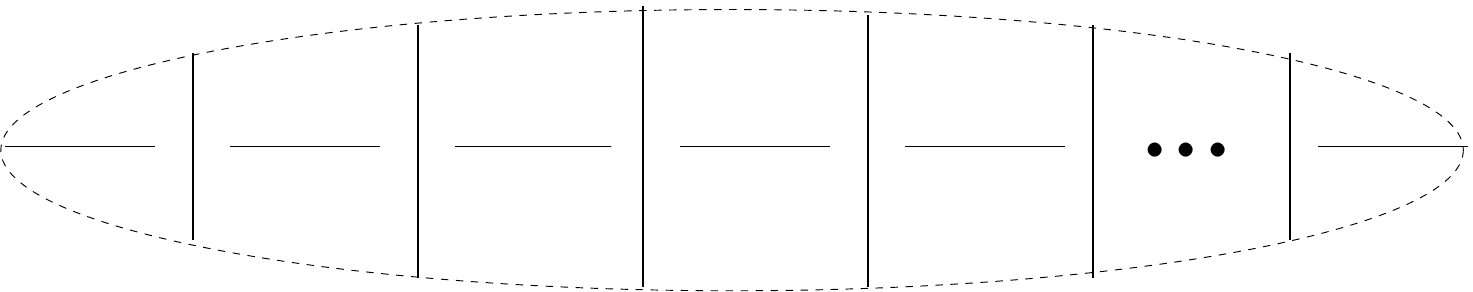}
\end{center}
\caption{A neighborhood of $U_1$}
\label{fig:ubridge1}
\end{figure}

\begin{figure}[ht!]
\begin{center}
\labellist
\pinlabel $a_0$ [r] at 452 257
\pinlabel $a_3$ [t] at 757 231
\pinlabel $b_0$ [l] at 306 253
\pinlabel $b_3$ [b] at 0 288
\pinlabel $U_1$ [r] at 505 299
\pinlabel $U_2$ [l] at 262 221
\pinlabel $O_1$ [bl] at 603 522
\pinlabel* $O_2$ [tr] at 73 60
\endlabellist
\includegraphics[width=3in]{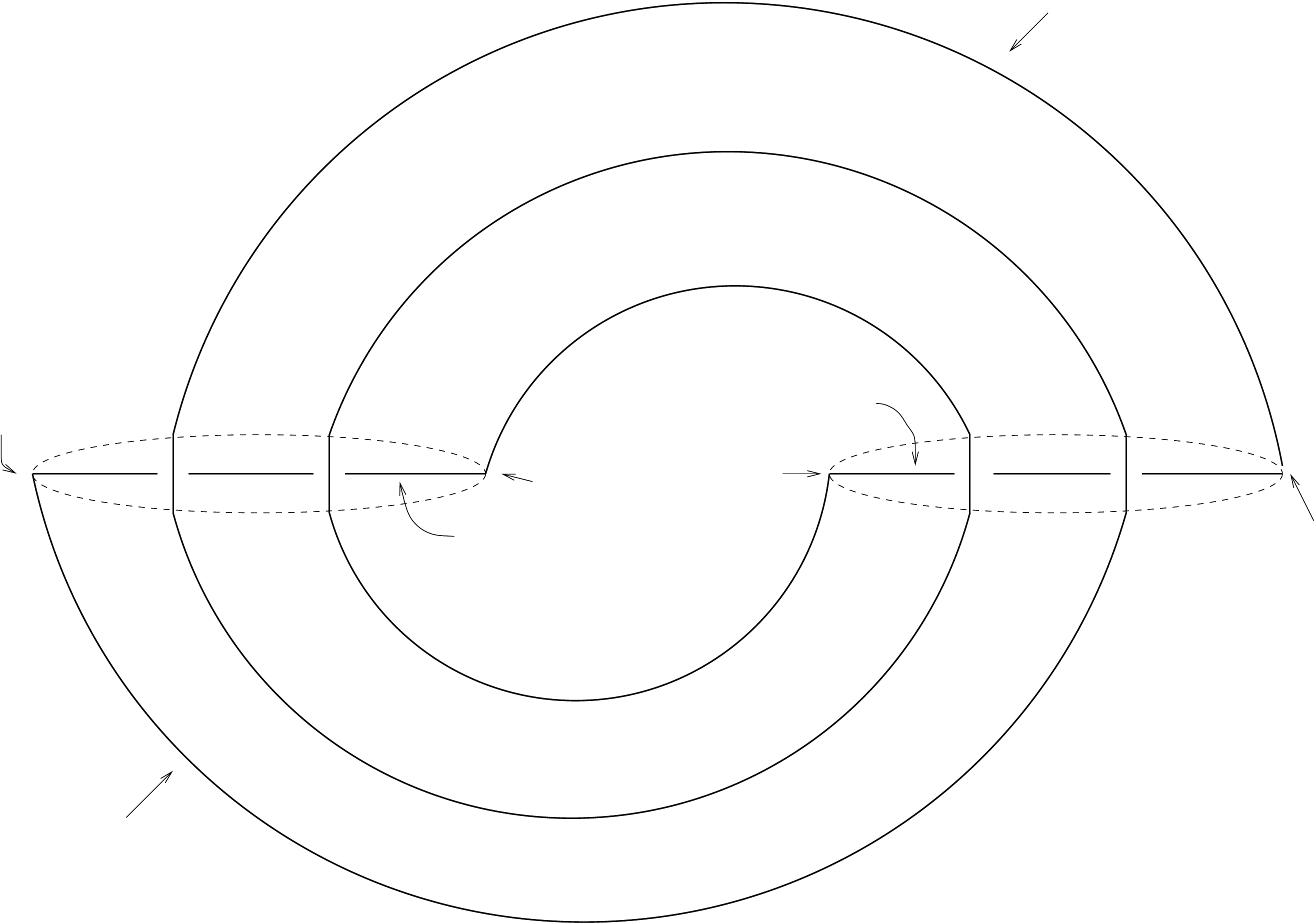}
\end{center}
\caption{Schubert normal form ($U_1 \cup U_2 \cup O_1 \cup O_2$) for 
$K(3,1)=$ the right-handed trefoil} 
\label{fig:obridge}
\end{figure}

  \item $O_1$ and $O_2$ are formed by connecting $a_{i \mod 2p}$ to $b_{(i-q) \mod 2p}$.  See \fullref{fig:obridge} for the example of $K(3,1)$.

Using Proposition 2.2 in \cite{MR1928176}, we get a genus $2$ handlebody decomposition for $K(p,q)$ that extends to a Heegaard decomposition of $S^3$.  Schematically, we can draw this handlebody decomposition for $S^3 - K$ by

\begin{enumerate}
  \item
    placing the feet of two $1$--handles at $a_0, a_p$ and $b_0, b_p$, respectively,
  \item
    letting $\alpha_1 = U_1 + $ core of the $1$ handle whose feet are at $a_0$ and $a_p$, pushed out to the Heegaard surface,
  \item 
    similarly letting $\alpha_2 = U_2 +$ core of the $1$--handle whose feet are at $b_0$ and $b_p$, pushed out to the Heegaard surface,
  \item
    letting $\beta_1 =$ boundary of a regular neighborhood of either $O_1$ or $O_2$.
\end{enumerate}

\begin{figure}[ht!]
\begin{center}
\labellist
\pinlabel $\alpha_2{=}U_2$ [b] at 98 340
\pinlabel $\alpha_1{=}U_1$ [t] <4pt,0pt> at 724 181
\pinlabel $\beta_1{=}\partial(N(O_1))$ [b] at 466 533
\endlabellist
\includegraphics[width=3in]{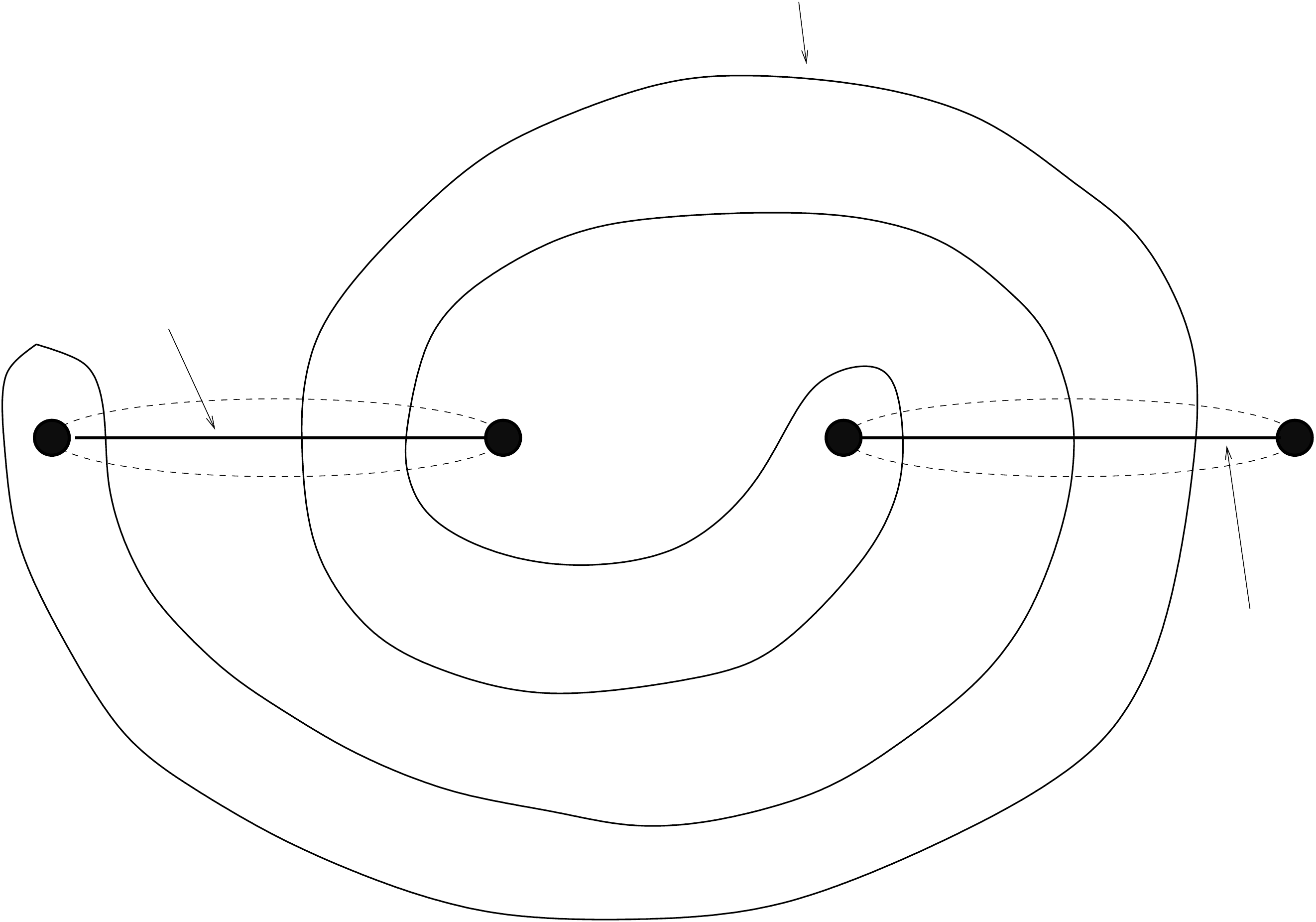}
\end{center}
\caption{A handlebody decomposition for $S^3 - K(3,1)$}
\label{fig:S3minusK}
\end{figure}

\fullref{fig:S3minusK} provides an illustration of this for $K=$ right-handed trefoil.

Hence, Schubert normal form for a two-bridge knot $K$ yields a genus $2$ doubly-pointed Heegaard diagram $\hd(S^3;K;\mu=\beta_2=h_{\alpha_2})$.  One easily checks that $h_{\alpha_1}, h_{\alpha_2}$ are both primitive in $H_1(S^3-K)$.

As described in \fullref{section:Handlebody}, this handlebody decomposition lifts to one for \sloppy \mbox{$\Sigma^2(K) - \widetilde{K}$} \fussy with an action of $\mathbb{Z}_2$.  If we denote the non-trivial element of $\mathbb{Z}_2$ by $\tau_2$, then $\hb(\Sigma^2(K)$ $- \widetilde{K}; \widetilde{\mu})$ has 

\begin{itemize}
  \item one $0$--handle $\widetilde{h}_0$
  \item three $1$--handles $\widetilde{h}_{\alpha_1}, \tau_2(\widetilde{h}_{\alpha_1}), \widetilde{h}_{\alpha_2}$ 
  \item two $2$--handles $\widetilde{h}_{\beta_1}, \tau_2(\widetilde{h}_{\beta_1})$.  
\end{itemize}

$\hb(\Sigma^2(K); \widetilde{K}; \widetilde{\mu} = \widetilde{h}_{\alpha_2})$ and $\hd(\Sigma^2(K); \widetilde{K})$ are as described in \fullref{section:Handlebody}.

We begin our $\widehat{HFK}(\Sigma^2(K); \widetilde{K})$ calculation by splitting generators into Spin$^c$ classes.  The following (very easy) observation will help us:

\begin{lemma} \label{lemma:transfer}
  Suppose $h$ is a $1$--cycle in $S^3$.  Let $\widetilde{h}$ and $\tau_2(\widetilde{h})$ be its two lifts in $\Sigma^2(K)$.  Then $\widetilde{h} + \tau_2(\widetilde{h}) = 0$ in $H_1(\Sigma^2(K))$.   
\end{lemma}

In particular, if $\widetilde{h}$ and $\tau_2(\widetilde{h})$ are themselves $1$--cycles in $\Sigma^2(K)$, then $\widetilde{h} = -\tau_2(\widetilde{h})$ in $H_1(\Sigma^2(K))$.

\proof[Proof of \fullref{lemma:transfer}]\qua
Let $$p_!(h) = \sum_{g \in \mathbb{Z}_2} g(\widetilde{h})$$ ($\widetilde{h}$ is a choice of lift of $h$) be the transfer map on chains (see, eg, Defn. 11.2 in \cite{MR1700700}) associated to the map $$p\co\Sigma^2(K) \rightarrow S^3.$$  Since $h = 0 \in H_1(S^3)$ and the transfer map is a homomorphism on homology, $p_!(h) = \widetilde{h} + \tau_2(\widetilde{h}) = 0 \in H_1(\Sigma^2(K))$, as desired. \endproof

\begin{theorem} \label{theorem:central}
  Given any knot $K$ in $S^3$ and a particular doubly-pointed Heegaard diagram $\hd(S^3;K)$, we can construct $\hd(\Sigma^2(K);\widetilde{K})$.  Then with respect to these particular Heegaard diagrams there is a natural map on chains $$f\co \widehat{CFK}(S^3;K) \rightarrow \widehat{CFK}(\Sigma^2(K);\widetilde{K})$$ given by $$f(\begin{bf}x \end{bf}) = (\widetilde{\begin{bf}x \end{bf}}, \tau_2(\widetilde{\begin{bf} x \end{bf}}))$$ for $\begin{bf} x \end{bf} \in \widehat{CFK}(S^3;K)$ .

All generators in the image of $f$ lie in the same Spin$^c$ structure, which we will denote $\mathfrak{s}_0$.

If $K$ is a two-bridge knot and $\hd(S^3;K)$ is the Heegaard diagram associated to the Schubert normal form of $K$, then $f$ is a chain map and induces an isomorphism on homology; ie, $$f_*\co \widehat{HFK}(S^3; K) \rightarrow \widehat{HFK}(\Sigma^2(K); \widetilde{K}, \mathfrak{s}_0)$$ is an isomorphism.  

\end{theorem}

Here ${\bf x}$ is a $g$--tuple of intersection points of $\alpha$ and $\beta$ curves in $\hd(S^3;K)$, $\widetilde{{\bf x}}$ is a choice of lift of that $g$--tuple in $\hd(\Sigma^2(K);\widetilde{K})$, and $\tau_2(\widetilde{{\bf x}})$ is the image of that $g$--tuple under the non-trivial deck transformation.

\medskip
{\bf Remark}\qua It is worthwhile to mention that in the case where $K$ is a two-bridge knot in $S^3$, $\mathfrak{s}_0$ is the unique spin element of Spin$^c(\Sigma^2(K))$.

More precisely, consider the first Chern class map $c_1\co$ Spin$^c(Y) \rightarrow H^2(Y;\mathbb{Z})$: $$c_1(\mathfrak{s}) = \mathfrak{s} - \bar{\mathfrak{s}},$$ where $\bar{\mathfrak{s}}$ is the conjugate Spin$^c$ structure (see \cite{MR2113019}, Section 2.6).  In our situation, the map $\tau_2$ is the conjugation map on Spin$^c$ structures: $\tau_2(\mathfrak{s}) = \bar{\mathfrak{s}}$ for all $\mathfrak{s} \in$ Spin$^c(Y)$.

But it is clear that $\mathfrak{s}_0 = \tau_2(\mathfrak{s}_0)$, so $$c_1(\mathfrak{s}_0)= c_1(\bar{\mathfrak{s}_0});$$ ie, $\mathfrak{s}_0$ is spin.  

Furthermore, $\mathfrak{s}_0$ is the unique Spin$^c$ structure with this property, for if $\mathfrak{s} \neq \mathfrak{s}_0$, we have $\mathfrak{s} \neq \tau_2(\mathfrak{s})$, implying $c_1(\mathfrak{s}) \neq 0$.  But $H^2(Y;\mathbb{Z}) \cong \mathbb{Z}_k$ ($k$ odd) has no $2$--torsion, so $$c_1(\mathfrak{s}) - c_1(\tau_2(\mathfrak{s})) = 2c_1(\mathfrak{s}) \neq 0,$$ so $c_1(\mathfrak{s}) \neq c_1(\bar{\mathfrak{s}})$.

\proof[Proof of \fullref{theorem:central}]\qua We begin by showing
that all generators of the form $f(\begin{bf}x \end{bf}) =
(\widetilde{\begin{bf}x \end{bf}}, \tau_2(\widetilde{\begin{bf} x
\end{bf}}))$ in $\widehat{CFK}(\Sigma^2(K);\widetilde{K})$ are in the
same Spin$^c$ structure.  Recall \cite{MR2113019} that two
generators ${\bf x}$ and ${\bf y}$ of $\widehat{HFK}(Y)$ lie in the
same Spin$^c$ structure iff there exists some path ${\bf x}
\rightarrow {\bf y}$ along $\alpha$ curves and some path ${\bf y}
\rightarrow {\bf x}$ along $\beta$ curves such that the union
represents the $0$ element in $H_1(Y)$.

In our situation, we have two generators, $(\widetilde{{\bf x}},\tau_2(\widetilde{{\bf x}}))$ and $(\widetilde{{\bf y}}, \tau_2(\widetilde{{\bf y}}))$, and we wish to show that we can find a $1$--cycle as above representing $0 \in H_1(\Sigma^2(K))$.  

But now note that for any two generators ${\bf x}$ and ${\bf y}$ of $\widehat{HFK}(S^3)$ we can find a path $\gamma_{\alpha}$ traveling from ${\bf x}$ to ${\bf y}$ along $\alpha$ curves and a path $\gamma_\beta$ traveling from ${\bf y}$ to ${\bf x}$ along $\beta$ curves in the Heegaard diagram for $S^3$.  The union, $\gamma = \gamma_\alpha \bigcup \gamma_\beta$, is a closed $1$--cycle representing the trivial (only) element in $H_1(S^3)$.

But $p_!(\gamma)$ (where, again, $p_{!}$ is the transfer map) in the
Heegaard diagram $\hd(\Sigma^2(K);$ $\widetilde{K}; \alpha_2)$ exactly
gives a path from $(\widetilde{{\bf x}}, \tau_2(\widetilde{{\bf x}}))$
to $(\widetilde{{\bf y}}, \tau_2(\widetilde{{\bf y}}))$ along $\alpha$
curves and a path from $(\widetilde{{\bf y}}, \tau_2(\widetilde{{\bf
y}}))$ to $(\widetilde{{\bf x}}, \tau_2(\widetilde{{\bf x}}))$ along
$\beta$ curves.

Since $p_!(\gamma) = 0$ in $H_1(\Sigma^2(K))$, all generators of the form $(\widetilde{{\bf x}}, \tau_2(\widetilde{{\bf x}}))$ lie in the same Spin$^c$ structure, which we have called $\mathfrak{s}_0$.

We now turn to showing that, in the case of a two-bridge knot, $$\widehat{HFK}(\Sigma^2(K); \widetilde{K}, \mathfrak{s}_0) \cong \widehat{HFK}(S^3; K).$$
We will do so by showing that the lifted Heegaard diagram associated to the Schubert normal form for K, $\hd(\Sigma^2(K); \widetilde{K}; \widetilde{h}_{\alpha_2})$, has the property that no other generators lie in $\mathfrak{s}_0$.  The differentials and filtrations in this central Spin$^c$ structure will match the differentials and filtrations downstairs. 

\begin{proposition} \label{proposition:bijection}
  Let $K$ be a two-bridge knot, $\hd(S^3;K,h_{\alpha_2})$ be the genus $2$ Heegaard diagram for $K$ obtained from Schubert normal form, and \sloppy $\hd(\Sigma^2(K);\widetilde{K},\widetilde{h}_{\alpha_2})$ be the lifted genus $3$ Heegaard diagram for $\widetilde{K}$ in $\Sigma^2(K)$.

Then the map $f\co \widehat{CFK}(S^3;K) \rightarrow \widehat{CFK}(\Sigma^2(K); \widetilde{K}; \mathfrak{s}_0)$ described in \fullref{theorem:central} is a bijection of sets.  In particular, all elements of $\widehat{CFK}(\Sigma^2(K);\mathfrak{s}_0)$ are of the form $(\widetilde{{\bf x}}, \tau_2(\widetilde{{\bf x}}))$ for $\begin{bf} x \end{bf} \in \widehat{CFK}(S^3; K)$.
\end{proposition}

\proof[Proof of \fullref{proposition:bijection}]\qua
We have already shown that all of the $\tau_2$--invariant generators (those of the form $(\begin{bf}\widetilde{x} \end{bf}, \tau_2(\widetilde{{\bf x}}))$) are in $\mathfrak{s}_0$.  We now need only show that no non-$\tau_2$--invariant generators are in $\mathfrak{s}_0$.

The proof will rely on the fact that for a two-bridge knot, $$rk(\widehat{CFK}(S^3;K)) = |H_1(\Sigma^2(K))|,$$ which will imply that no non-$\tau_2$--invariant generators can appear in $\mathfrak{s}_0$.\footnote{This condition on $rk(\widehat{CFK}(S^3;K))$ holds for a wider class of knots (eg, alternating knots), and therefore, similar results may hold in wider generality.}
 
First, note that for $K$ a two-bridge knot with the standard genus $2$ Heegaard decomposition given by Schubert normal form, all generators in $\widehat{CFK}(\Sigma^2(K);\widetilde{K})$ for the genus $3$ Heegaard diagram constructed as the lift of this genus $2$ Heegaard diagram (as described in \fullref{subsection:HBody}) are naturally of the form $(\widetilde{{\bf x}}, \tau_2(\widetilde{{\bf y}}))$, where ${\bf x}$ and ${\bf y}$ are generators in $\widehat{CFK}(S^3)$, and $\widetilde{{\bf x}}, \tau_2(\widetilde{{\bf y}}) \in \mathbb{T}_{\widetilde{\begin{bf}\alpha \end{bf}}} \cap \mathbb{T}_{\widetilde{\begin{bf} \beta \end{bf}}}$ are lifts of ${\bf x}$ and ${\bf y}$.

More precisely, note that the generators in $\widehat{CFK}(S^3; K)$ are naturally identified with intersection points of $\alpha_1$ with $\beta_1$ (since $\beta_2$, the chosen meridian of the knot, intersects only $\alpha_2$).  Similarly,  the generators of $\, \widehat{CFK}(\Sigma^2(K); \widetilde{K})$ in the lifted Heegaard diagram for $\Sigma^2(K)$, are in one-to-one correspondence with formal summands of the determinant 
$$\begin{array}{ll}
& \begin{array}{rl} \hskip 8pt \widetilde{\beta}_1 & \hskip 5pt \tau_2(\widetilde{\beta}_1) \end{array} \\
\begin{array}{r} \hskip 2pt \widetilde{\alpha}_1 \\ \tau_2(\widetilde{\alpha}_1) \end{array} & \left( \begin{array}{cr} &  \\  &  \hskip 40pt \end{array} \right) \end{array}
$$

where the entries in the matrix above are formal sums of intersection points of the corresponding $\alpha$ and $\beta$ curves.  A summand of this matrix can be thought of as a pair $(\widetilde{{\bf x}}, \tau_2(\widetilde{{\bf y}}))$ where ${\bf x}$, ${\bf y}$ are generators of $\widehat{CFK}(S^3;K)$ and $\widetilde{{\bf x}}$ and $\widetilde{{\bf y}}$ are the lifts which lie on the $\widetilde{\alpha}_1$ curve.\footnote{Such an identification of generators of $\widehat{CFK}(\Sigma^2(K);\widetilde{K})$ with pairs of generators of $\widehat{CFK}(S^3;K)$ is possible for a general choice of Heegaard diagram for $K$ in $S^3$, but different generators will require different lifts of $\mathbb{T}_{\alpha}$ and $\mathbb{T}_{\beta}$ in order to make the identification.}

We can then measure the Spin$^c$ structure of a non-$\tau_2$--invariant generator $(\widetilde{\begin{bf}x \end{bf}}, \tau_2(\widetilde{{\bf y}}))$ by comparing it to the $\tau_2$--invariant generator $(\widetilde{\begin{bf} y \end{bf}}, \tau_2(\widetilde{{\bf y}}))$.  More precisely:

\begin{lemma} \label{lemma:Spinc=pi1}
If $(\widetilde{{\bf x}}, \tau_2(\widetilde{{\bf y}}))$ is a generator of $\widehat{CFK}(\Sigma^2(K),\widetilde{K})$, then its corresponding element of Spin$^c(\Sigma^2(K))$ relative to $\mathfrak{s}_0$, ie, $$\mathfrak{s}(\widetilde{{\bf x}}, \tau_2(\widetilde{{\bf y}})) - \mathfrak{s}(\widetilde{{\bf y}}, \tau_2(\widetilde{{\bf y}})),$$ thought of as an element of $H_1(\Sigma^2(K)),$ is represented by the lift, $\widetilde{\gamma}$, of the word $\gamma$ in $\pi_1(S^3-K)$ read off as we travel from ${\bf x}$ to ${\bf y}$ along $\beta_1.$
\end{lemma}

\proof[Proof of \fullref{lemma:Spinc=pi1}]\qua
$\mathfrak{s}(\widetilde{{\bf x}},\tau_2(\widetilde{{\bf y}})) - \mathfrak{s}(\widetilde{{\bf y}},\tau_2(\widetilde{{\bf y}}))$, as an element of $H_1(\Sigma^2(K))$, is represented by the cycle obtained by connecting $(\widetilde{{\bf x}},\tau_2(\widetilde{{\bf y}}))$ to $(\widetilde{{\bf y}},\tau_2(\widetilde{{\bf y}}))$ along $\alpha$ curves and $(\widetilde{{\bf y}},\tau_2(\widetilde{{\bf y}}))$ to $(\widetilde{{\bf x}},\tau_2(\widetilde{{\bf y}}))$ along $\beta$ curves (see Definition 2.4 and Section 2.6 of \cite{MR2113019}; also see \cite{MR1484699}).

We construct such a path as the product of:
\begin{itemize}
\item
the appropriate lift of a loop between ${\bf x}$ and ${\bf y}$ (along $\alpha_1$ and back along $\beta_1$) in $\hd(S^3;K)$ to a loop between $\widetilde{{\bf x}}$ and $\widetilde{{\bf y}}$ in $\hd(\Sigma^2(K);\widetilde{K})$,
\item
the constant path from $\tau_2(\widetilde{{\bf y}})$ to $\tau_2(\widetilde{{\bf y}})$,
\item
a path from the lone intersection point between $\widetilde{\beta}_2$ and $\widetilde{\alpha}_2$ to itself along $\widetilde{\beta}_2$.
\end{itemize}

The element of $H_1(\Sigma^2(K))$ represented by this path is represented by the word in $\pi_1(\Sigma^2(K))$ read off as we travel from $\widetilde{{\bf x}}$ to $\widetilde{{\bf y}}$ along $\widetilde{\beta_1}$, which is the lift, $\widetilde{\gamma}$, of the word $\gamma$ in $\pi_1(S^3-K)$ read off as we travel from ${\bf x}$ to ${\bf y}$ along $\beta_1$. \endproof

But now I claim that, since $rk(\widehat{CFK}(S^3;K)) = |H_1(\Sigma^2(K))|$ for the Heegaard diagram associated to Schubert normal form for $K$ a two-bridge knot (see \cite{GT0306378}, \cite{MR1928176}, \cite{MR1988285}), $\widetilde{\gamma} \neq 0$ unless ${\bf x} = {\bf y}$.

\begin{lemma} \label{lemma:Rk=Det}
Let $K$ a two-bridge knot in $S^3$ and  ${\bf x}, {\bf y}$ two generators of $\widehat{CFK}(S^3;K)$ associated to the handlebody decomposition coming from Schubert normal form.  Let $\gamma({\bf x},{\bf y})$ be the word in $\pi_1(S^3-K)$ read off as we travel from ${\bf x}$ to ${\bf y}$ along $\beta_1$ and $\widetilde{\gamma}({\bf x},{\bf y})$ a lift of $\gamma$ to a word in $\pi_1(\Sigma^2(K)-\widetilde{K})$.

Then ${\bf x} \neq {\bf y}$ implies that $\widetilde{\gamma}({\bf x},{\bf y}) \neq 0$ as an element of $H_1(\Sigma^2(K))$.
\end{lemma}

\proof[Proof of \fullref{lemma:Rk=Det}]\qua
The crucial observation is that the Fox matrix associated to the homomorphism $$\pi_1(S^3-K) \rightarrow \mathbb{Z}_2$$ is a presentation matrix for $H_1(\Sigma^2(K))$, and its summands are in natural one-to-one correspondence with the generators of $\widehat{CFK}(S^3;K)$. See \fullref{subsection:determinant} for a more detailed discussion of Fox calculus.

Fix a generator, ${\bf x}$, of $\widehat{CFK}(S^3;K)$ and begin reading off the relation corresponding to the boundary of $\beta_1$, beginning at ${\bf x}$.  Let $$t\co \pi_1(S^3-K) \rightarrow \mathbb{Z}_2$$ be the homomorphism inducing the branched double cover of $K$ and $$\tau\co \mathbb{Z}[\pi_1(S^3-K)] \rightarrow \mathbb{Z}[\mathbb{Z}_2]$$ the group-ring extension.

Recall that the generators of $\widehat{CFK}(S^3;K)$ correspond one-to-one with the intersection points of $\alpha_1$ and $\beta_1$.  For convenience, label the generators of $\widehat{CFK}(S^3;K)$ by ${\bf y}_i$, according to the order in which we encounter them as we travel along $\beta_1$ from ${\bf x}$.

Then $$\frac{\partial \beta_1}{\partial \alpha_1} = \sum_{{\bf y}_i \in \widehat{CFK}(S^3;K)}\pm \gamma({\bf x}, {\bf y}_i);$$ ie, the Fox derivative of $\beta_1$ by $\alpha_1$ is the formal sum of the words connecting ${\bf x}$ to each of the other generators of $\widehat{CFK}(S^3;K)$, with signs in the sum given by the local intersection number of $\alpha_1 \cap \beta_1$ at ${\bf y}_i$.  Notice that there are exactly $rk(\widehat{CFK}(S^3;K))$ summands in the Fox determinant. 

Furthermore, $\tau(\frac{\partial \beta_1}{\partial \alpha_1})$ is a presentation matrix for $H_1(\Sigma^2(K))$, after replacing the formal elements of $\mathbb{Z}_2$ by square roots of unity ($\pm 1$).

But $rk(\widehat{CFK}(S^3;K)) = |H_1(\Sigma^2(K))|$, so all summands must have the same sign once we replace the formal elements of $\mathbb{Z}_2$ with $\pm 1$ and take into account local intersection multiplicities.

Since the element of $H_1(\Sigma^2(K))$ represented by $\widetilde{\gamma}({\bf x},{\bf y}_k)$ is precisely $$\sum_{i=1}^k \tau(\pm \gamma({\bf x},{\bf y}_i)) \cdot h_{\widetilde{\alpha}_1},$$ and $h_{\widetilde{\alpha}_1}$ is a primitive generator of $H_1(\Sigma^2(K))$, we conclude that $\widetilde{\gamma}({\bf x},{\bf y})$ cannot be $0$ in $H_1(\Sigma^2(K))$ for any ${\bf y} \neq {\bf x}$.
\endproof

The following two lemmas prove that the map $$f\co \widehat{CFK}(S^3;K) \rightarrow \widehat{CFK}(\Sigma^2(K); \widetilde{K}; \mathfrak{s}_0)$$ given above also preserves the relative filtration and homological gradings in the case that $K$ is a two-bridge knot.  Hence, $f$ is a chain map respecting the filtration, so $f_*\co \widehat{HFK}(S^3; K) \rightarrow \widehat{HFK}(\Sigma^2(K); \widetilde{K}; \mathfrak{s}_0)$ is an isomorphism. 

\begin{lemma} \label{lemma:sgrading}
  For $K$ any knot in $S^3$, $$s(\begin{bf} x \end{bf}, \begin{bf} y \end{bf}) = s(f(\begin{bf} x \end{bf}), f(\begin{bf} y \end{bf}))$$ for all pairs of generators ${\bf x}$ and ${\bf y}$ in $\widehat{HFK}(S^3; K)$.
\end{lemma}

\proof[Proof of \fullref{lemma:sgrading}]\qua
Let ${\bf x}$ and ${\bf y}$ be two generators in $\widehat{HFK}(S^3;K)$ and $\phi \in \pi_2({\bf x}, {\bf y})$ a topological disk with $n_z(\phi) - n_w(\phi) = 1$.  

We have already observed in the proof of \fullref{theorem:central} that if $\gamma$ is the boundary of the image of $\phi$ in the Heegaard surface, $S$, then $\widetilde{\gamma} = p_!(\gamma)$ will be the boundary of the image of $\widetilde{\phi} \in \pi_2(f(\begin{bf} x \end{bf}), f(\begin{bf} y \end{bf}))$ in $\widetilde{S}$.

\fullref{lemma:filtration} asserts that $s({\bf x}, {\bf y})$ is equal to the coefficient on $\mu$, our choice of meridian, in $\gamma$.

But by the way we have defined $\hd(\Sigma^2(K); \widetilde{K}, \widetilde{\mu})$, the coefficient on $\widetilde{\mu}$ in $\widetilde{\gamma}$ is equal to the coefficient on $\mu$ in $\gamma$.  \endproof

\begin{lemma} \label{lemma:mgrading}
  For $K$ a two-bridge knot in $S^3$ and $\hd(S^3;K)$ a Heegaard diagram for $K$ coming from Schubert normal form as before, $$m(\begin{bf} x \end{bf}, \begin{bf}y \end{bf}) = m(f(\begin{bf} x \end{bf}), f(\begin{bf} y \end{bf})).$$
\end{lemma}

\proof[Proof of \fullref{lemma:mgrading}]
Consider $\hd(S^3;K)$ and $\hd(\Sigma^2(K); \widetilde{K}, \widetilde{h}_{\alpha_2})$ obtained from Schubert normal form.  It will be convenient to destabilize each of these Heegaard diagrams once by canceling the $h_{\alpha_2}, h_{\beta_2}$ pair and the $\widetilde{h}_{\alpha_2}, \widetilde{h}_{\beta_2}$ pair, respectively.  Call these destabilized Heegaard diagrams $\hd^\circ(S^3;K)$ and $\hd^\circ(\Sigma^2(K); \widetilde{K}, \widetilde{h}_{\alpha_2})$, respectively.

Note that there is a natural one-to-one correspondence between $\widehat{CFK}$ generators corresponding to $\hd$ and $\hd^\circ$ (see Prop. 6.1 in \cite{MR2065507}) which induces the isomorphism on $\widehat{HFK}$ corresponding to the (de)stabilization.

This destabilized Heegaard diagram for $K$ coming from Schubert normal form is particularly nice because we can find relative Maslov gradings between all generators just by looking at disks coming from ``finger moves.''

Namely, we know that we obtain $\hd^\circ(S^3; K)$ from the standard genus $1$ Heegaard diagram for $S^3$ by performing finger moves of the curve $\beta$ across the curve $\alpha$.

In fact, one can check that for each generator, there is a natural disk connecting it to at least one adjacent generator by a finger.  Each of these disks has either $n_w = 1, n_z = 0$ or $n_w = 0, n_z = 1$, and these ``finger disks'' are enough to determine the relative Maslov grading of any two generators.

Now, suppose that ${\bf x}$ and ${\bf y}$ are two generators in $\widehat{CFK}(S^3;K)$ connected by a finger disk.  Then $(\widetilde{\begin{bf} x \end{bf}}, \tau_2(\widetilde{\begin{bf} x \end{bf}}))$ and $(\widetilde{\begin{bf}y \end{bf}}, \tau_2(\widetilde{\begin{bf}y \end{bf}}))$ are connected by the lift of the finger disk, which is a quadrilateral.  Such a quadrilateral always represents a holomorphic disk in $\Sym^2(S)$ of Maslov index $1$.

Therefore, $$m(f(\begin{bf} x \end{bf}), f(\begin{bf} y \end{bf})) = m({\bf x}, {\bf y})$$ whenever ${\bf x}$ and ${\bf y}$ are connected by a finger disk.

Now we claim that every generator is connected to every other generator by some series of finger disks.  We can see this as follows:

Notice that $\beta_1$ in $\hd^{\circ}(S^3;K)$ is the image, after destabilizing the $\hd(S^3;K)$ coming from Schubert normal form, of a regular neighborhood of either one of the overbridges (note that a regular neighborhood of $O_1$ is isotopic to a regular neighborhood of $O_2$).  Similarly, $\alpha_1$ is the union of one of the underbridges (whichever one was not canceled in the destabilization) -- $U_1$, say -- with the core of the $1$--handle whose feet are placed at $a_0$ and $a_p$.  

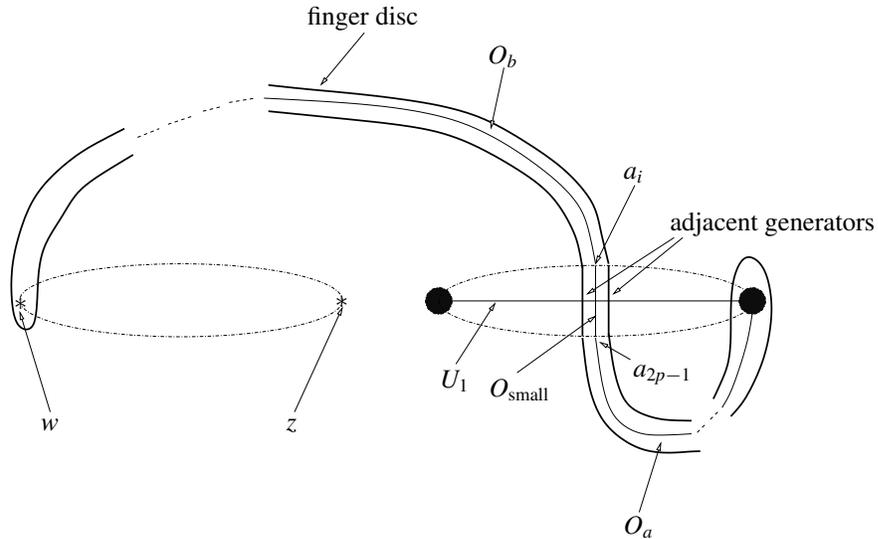
\begin{figure}[ht!]
\begin{center}
\labellist
\pinlabel $w$ [t] at 25 61
\pinlabel $z$ [t] at 175 61
\pinlabel {finger disc} [b] at 219 297
\pinlabel {adjacent generators} [bl] <-5pt,0pt> at 416 170
\pinlabel $a_{2p{-}1}$ [tl] at 385 94
\pinlabel $a_i$ [b] at 387 202
\pinlabel $*$ at 206 131
\pinlabel $*$ at 7 130
\pinlabel $U_1$ [t] at 276 91
\pinlabel $O_a$ [t] at 390 0
\pinlabel $O_b$ [b] at 306 274
\pinlabel $O_{\rm small}$ [t] <2pt,0pt> at 313 85
\endlabellist
\includegraphics[width=4in]{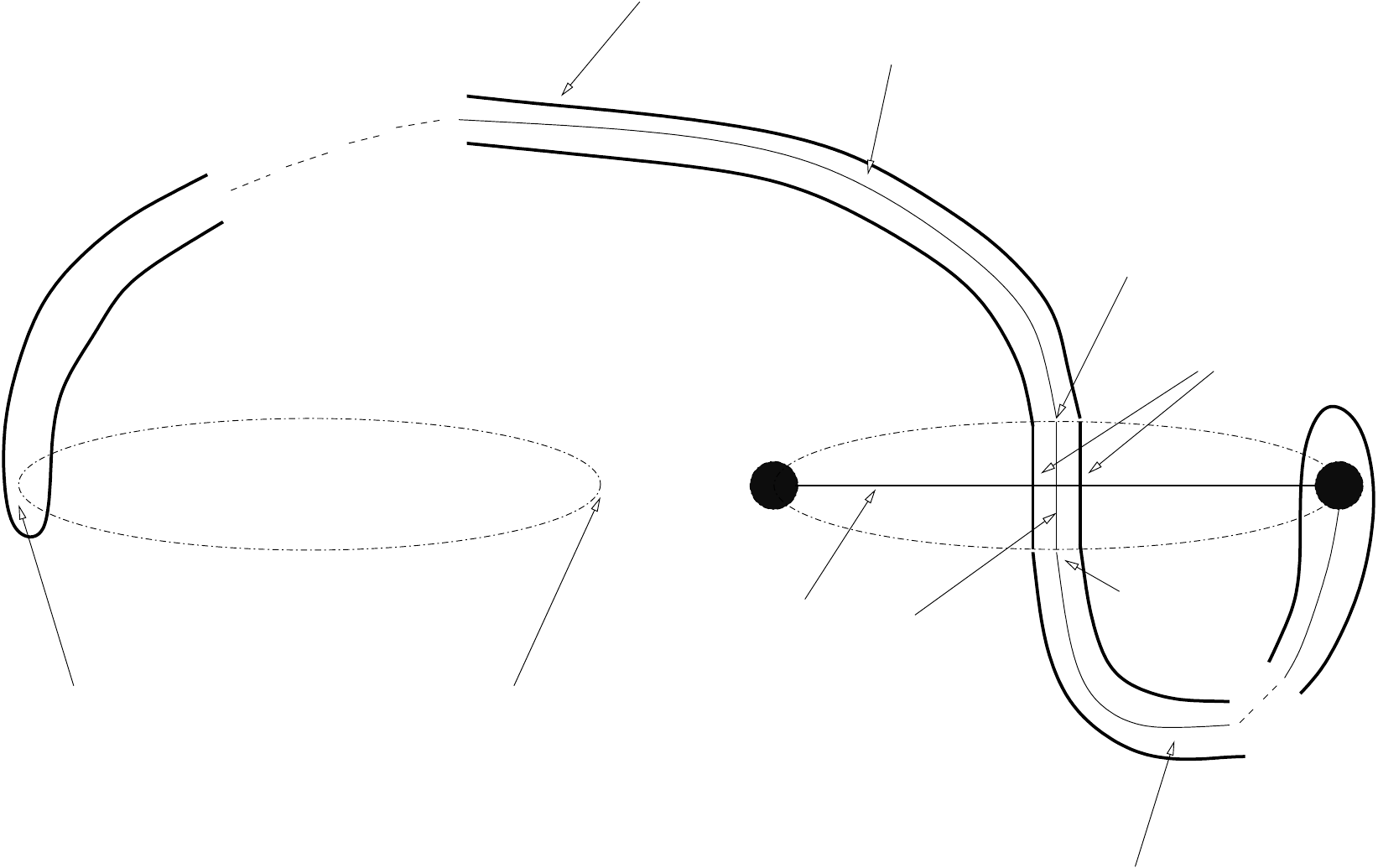}
\end{center}
\caption{Adjacent generators in a genus $1$ (destabilized) Heegaard diagram compatible with a two-bridge knot}
\label{fig:Fingerproof}
\end{figure}

Now, focus on any two adjacent generators of the Heegaard Floer homology (see \fullref{fig:Fingerproof}).  These two adjacent generators correspond to where $\beta_1$, considered either as $N(O_1)$ or $N(O_2)$, intersects $\alpha_1$.  For definiteness, say that the two adjacent generators we are considering are associated to $N(O_1)$.  Note that $O_1$ can be decomposed as the union of 

\begin{itemize} 
  \item a small overbridge $O_{small}$, which exits a neighborhood of $U_1$ at $a_i$ and $a_{2p - i}$ 
  \item an arc, $O_a$, connecting one of $a_i$ and $a_{2p-i}$ with an endpoint of $U_1$ (either $a_0$ or $a_p$), and 
  \item an arc, $O_b$, connecting the other endpoint ($a_{2p-i}$ or $a_i$) with an endpoint of $U_2$ (either $b_0$ or $b_p$).
\end{itemize}

Since $b_p$ and $b_0$ are where the basepoints $w$ and $z$ are positioned in the destabilized Heegaard diagram, a regular neighborhood of the arc connecting $a_i$ (or $a_{2p-i}$) to $b_0$ (or $b_p$) will be a finger disk connecting these two adjacent generators.

In other words, every two adjacent generators is connected by a finger disk.  So every pair of generators is connected by some sequence of finger disks.  

Therefore, $$m(f(\begin{bf} x \end{bf}), f(\begin{bf} y \end{bf}))= m({\bf x}, {\bf y})$$ for all pairs of generators ${\bf x}, {\bf y}$ in $\widehat{CFK}(S^3; K)$, as desired. \endproof

\subsection{Examples: $K(15,7)$ and $K(15,4)$} \label{section:examples}

\begin{figure}[ht!]
\begin{center}
\labellist
\pinlabel $x_1$ [t] at -26 68
\pinlabel $x_2$ [t] at 0 68
\pinlabel $x_3$ [t] at 26 68
\pinlabel $x_4$ [t] at 52 68
\pinlabel $x_5$ [b] at 303 265
\pinlabel $x_6$ [b] at 278 265
\pinlabel $x_7$ [b] at 253 265
\pinlabel $x_8$ [b] at 227 265
\pinlabel $y_1$ at 94 170
\pinlabel $y_7$ at 114 170
\pinlabel $y_2$ at 135 170
\pinlabel $y_6$ at 155 170
\pinlabel $y_3$ at 177 170
\pinlabel $y_5$ at 194 170
\pinlabel $y_4$ at 220 170
\pinlabel {filtration level} [t] at 14 0
\pinlabel ${-}1$ <0pt,-15pt> at 14 0
\pinlabel {filtration level} [t] at 157 0
\pinlabel $0$  <0pt,-15pt> at 157 0
\pinlabel {filtration level} [t] <10pt,0pt> at 267 0
\pinlabel $1$ <10pt,-15pt> at 267 0
\endlabellist
\includegraphics[width=2.5in]{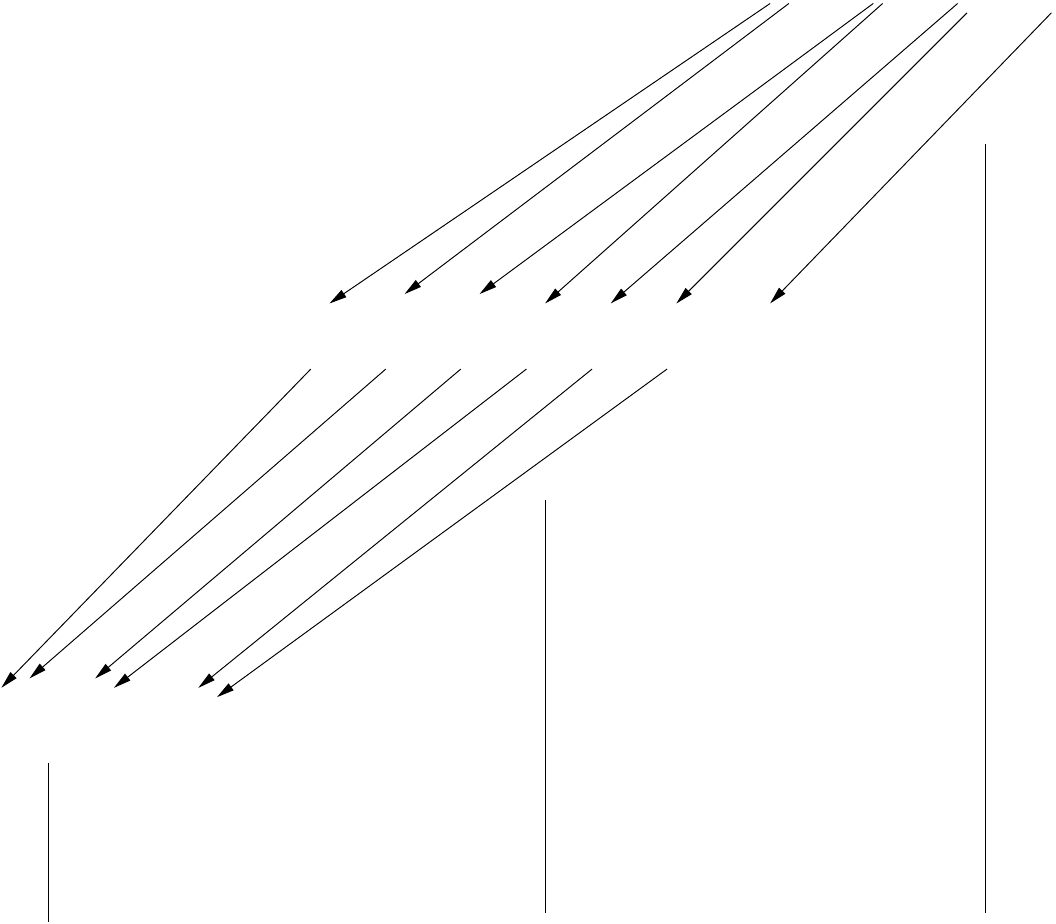}
\end{center}
\vspace{2mm}
\caption{$\mathbb{Z}$--filtered chain complex for $\widehat{CFK}(S^3; K(15,7))$}
\label{fig:K157cfhat}
\end{figure}

\begin{figure}[ht!]
\begin{center}
\labellist
\pinlabel $x_1$ [b] at 6 298
\pinlabel $x_2$ [b] at 64 298
\pinlabel $x_3$ [b] at 119 298
\pinlabel $x_4$ [b] at 177 298
\pinlabel $x_5$ [b] at 272 298
\pinlabel $x_6$ [b] at 335 298
\pinlabel $x_7$ [b] at 407 298
\pinlabel $x_8$ [b] at 469 298
\pinlabel $y_1$ [t] at 443 244
\pinlabel $y_7$ [t] at 33 244
\pinlabel $y_2$ [t] at 375 244
\pinlabel $y_6$ [t] at 90 244
\pinlabel $y_3$ [t] at 311 244
\pinlabel $y_5$ [t] at 145 244
\pinlabel $y_4$ [t] at 226 244
\pinlabel $w$ [t] at 222 293
\pinlabel $z$ [t] at 272 212
\pinlabel $*$ [b] at 227 315
\pinlabel $*$ [b] at 272 234
\endlabellist
\includegraphics[width=3.3in]{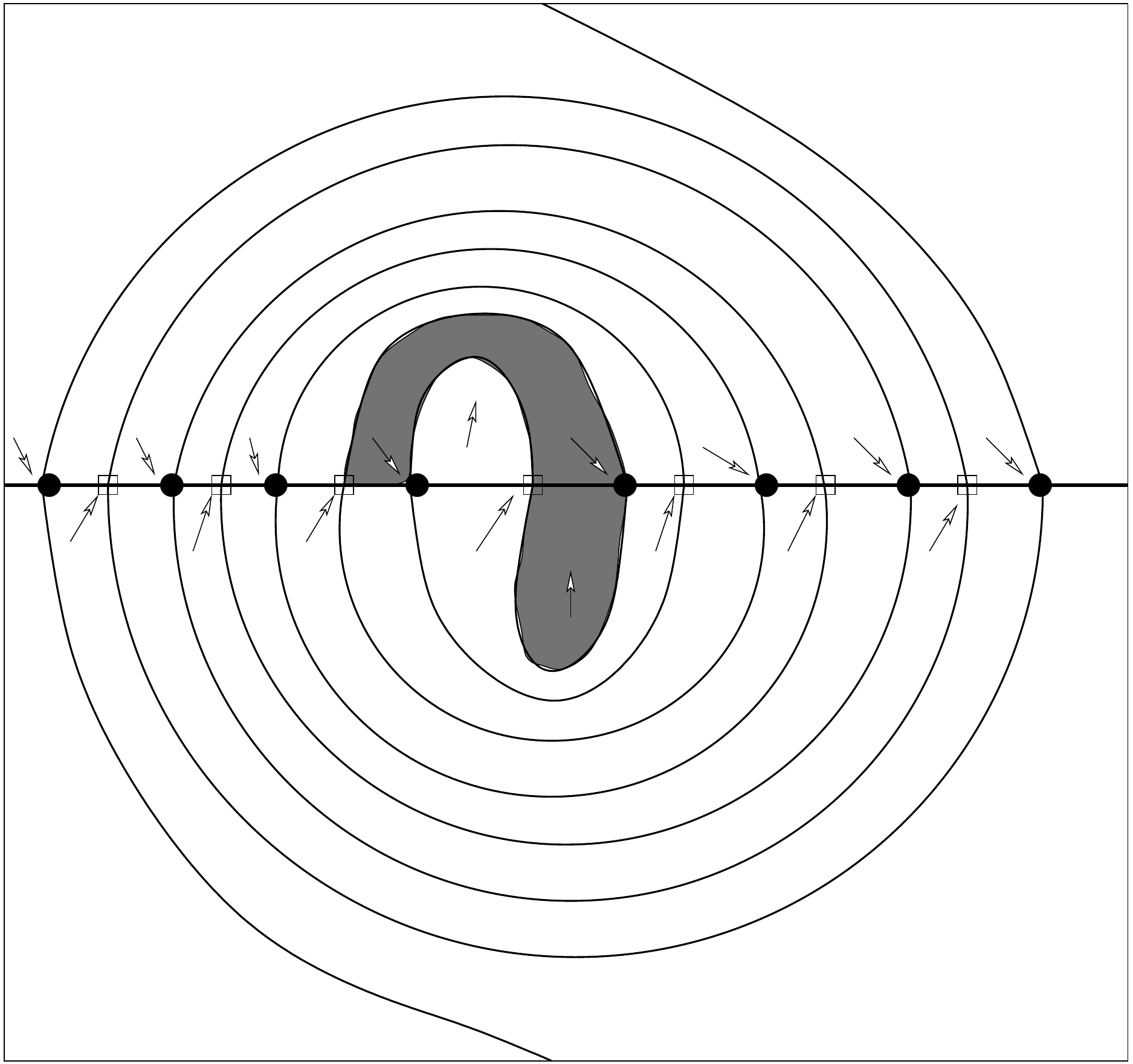}
\end{center}
\caption{Genus $1$ (destabilized) Heegaard diagram for $K(15,7)$ in $S^3$}
\label{fig:K157tinv}
\end{figure}

\begin{figure}[ht!]
\begin{center}
\labellist
\pinlabel $w$ [r] at 0 430
\pinlabel $\beta$ [r] at 171 301
\pinlabel $z$ [r] at 228 187
\pinlabel $x_1$ [b] <-2pt,0pt> at 288 83
\pinlabel $\tau_2(\beta)$ [b] at 451 361
\pinlabel $\tau_2(x_1)$ [t] <5pt,0pt> at 427 425
\pinlabel $\tau_2(\alpha)$ [l] at 500 545
\pinlabel $\alpha$ [l] at 497 18
\pinlabel $*$ at 46 403
\pinlabel $*$ at 212 569
\pinlabel $*$ at 447 569
\pinlabel $*$ at 613 402
\pinlabel $*$ at 613 167
\pinlabel $*$ at 445 1
\pinlabel $*$ at 211 1
\pinlabel $*$ at 46 168
\endlabellist
\includegraphics[width=3.4in]{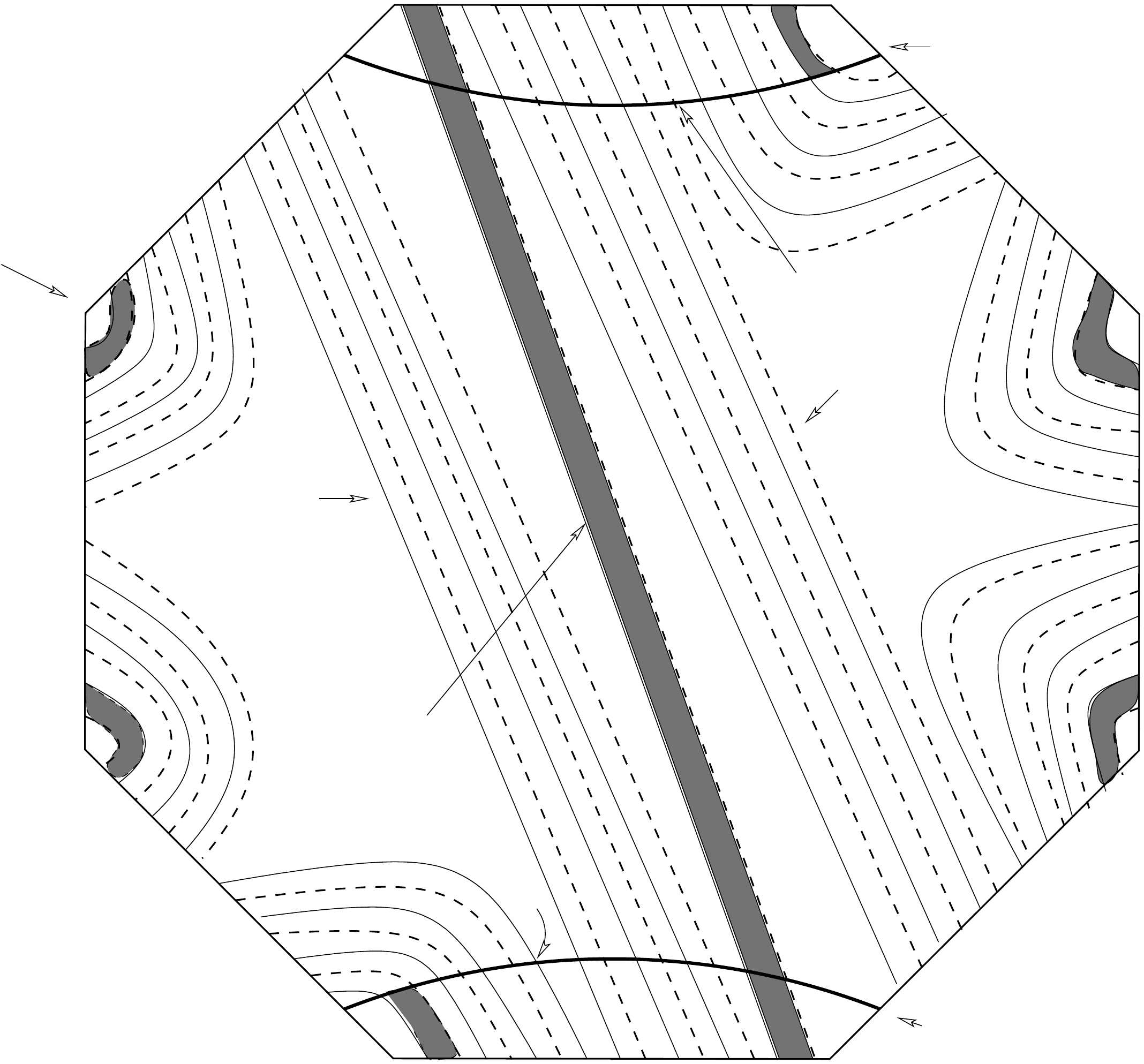}
\end{center}
\caption{Genus $2$ (destabilized) Heegaard diagram for $\widetilde{K}(15,7)$ in $L(15,7)$}
\label{fig:K157Branchedtinv}
\end{figure}

{\small
\begin{table}[ht!]
\centering
\vspace{1mm}
\begin{tabular}{|c|c|c|c|c|}
\hline
$\mathfrak{s}_0$ & $\mathfrak{s}_{\pm 1}$ & $\mathfrak{s}_{\pm 2}$ & $\mathfrak{s}_{\pm 3}$ & $\mathfrak{s}_{\pm 4}$\\
\hline
$(x_1, x_1)$, $(y_1, y_1)$& $(x_1,x_8)$ & $(x_1,x_2), (y_1,y_2)$ & $(x_1,x_7)$ &    $(x_1,x_3), (y_1,y_3)$\\ 
$(x_2, x_2)$, $(y_2, y_2)$&             & $(x_2,x_3), (y_2,y_3)$ & $(x_2,x_8)$ &    $(x_2,x_4), (y_2,y_4)$\\
$(x_3, x_3)$, $(y_3, y_3)$&             & $(x_3,x_4), (y_3,y_4)$ & $(y_1,y_7)$ &    $(x_3,x_5), (y_3,y_5)$\\
$(x_4, x_4)$, $(y_4, y_4)$&             & $(x_4,x_5), (y_4,y_5)$ &             &    $(x_4,x_6), (y_4,y_6)$\\
$(x_5, x_5)$, $(y_5, y_5)$&             & $(x_5,x_6), (y_5,y_6)$ &             &    $(x_5,x_7), (y_5,y_7)$ \\
$(x_6, x_6)$, $(y_6, y_6)$&             & $(x_6,x_7), (y_6,y_7)$ &             &    $(x_6,x_8)$           \\
$(x_7, x_7)$, $(y_7, y_7)$&             & $(x_7,x_8), (y_7,y_8)$ &             &                          \\         
$(x_8, x_8)$,             &             &                        &             &                          \\
\hline
\end{tabular}

\vspace{1mm}
\begin{tabular}{|c|c|c|}
\hline
$\mathfrak{s}_{\pm 5}$ & $\mathfrak{s}_{\pm 6}$ & $\mathfrak{s}_{\pm 7}$ \\
\hline
$(x_1,x_6)$ & $(x_1,x_4), (y_1,y_4)$ & $(x_1,x_5), (y_1,y_5)$\\ 
$(y_1,y_6)$ & $(x_2,x_5), (y_2,y_5)$ & $(x_2,x_6), (y_2,y_6)$\\
$(x_2,x_7)$ & $(x_3,x_6), (y_3,y_6)$ & $(x_3,x_7), (y_3,y_7)$\\
$(y_2,y_7)$ & $(x_4,x_7), (y_4,y_7)$ & $(x_4,x_8)$           \\
$(x_3,x_8)$ & $(x_5,x_8)$            &                       \\
\hline
\end{tabular}
\vspace{1mm}
\caption{Spin$^c$ structures $\mathfrak{s}_0, \ldots \mathfrak{s}_{\pm 7}$ for $\Sigma^2(K(15,7))$}
\label{table:Spinc157}
\end{table}
}

\begin{figure}[ht!]
\begin{center}
\labellist
\pinlabel $x_1$ [b] at 2 116
\pinlabel $x_2$ [b] at 25 116
\pinlabel $x_3$ [b] at 46 116
\pinlabel $x_4$ [b] at 69 116
\pinlabel $x_5$ [b] at 105 116
\pinlabel $x_6$ [b] at 131 116
\pinlabel $x_7$ [b] at 159 116
\pinlabel $x_8$ [b] at 183 116
\pinlabel $y_7$ [t] at 15 95
\pinlabel $y_6$ [t] at 37 95
\pinlabel $y_5$ [t] at 57 95
\pinlabel $y_4$ [t] at 88 95
\pinlabel $y_3$ [t] at 121 95
\pinlabel $y_2$ [t] at 146 95
\pinlabel $y_1$ [t] at 172 95
\pinlabel $w$ [t] at 86 114
\pinlabel $z$ [t] at 117 70
\pinlabel $*$ [b] at 88 123
\pinlabel $*$ [b] at 108 87
\endlabellist
\includegraphics[width=3.3in]{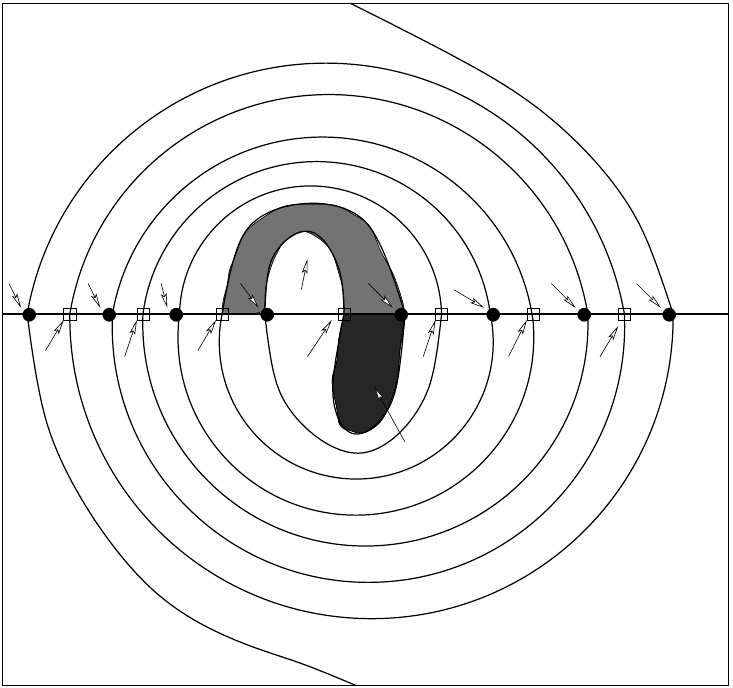}
\end{center}
\caption{Two holomorphic disks, one in $\pi_2(x_5,y_4)$ and the other in $\pi_2(y_5,x_4)$}
\label{fig:K157tnon}
\end{figure}

\begin{figure}[ht!]
\begin{center}
\labellist
\pinlabel $w$ [r] at 0 237
\pinlabel $*$ [lb] at 174 156
\pinlabel $\beta$ [r] at 93 168
\pinlabel $z$ [rt] at 126 99
\pinlabel $\tau_2(\beta)$ [b] at 247 199
\pinlabel $\tau_2(\alpha)$ [bl] at 265 314
\pinlabel $\alpha$ [l] at 273 12
\pinlabel {$\tilde x_4$} [r] at 71 37
\pinlabel {$\tilde y_5$} [r] at 97 0
\pinlabel {$\tilde y_4$} [b] at 242 32
\pinlabel {$\tilde x_5$} [b] at 228 60
\pinlabel {$\tau_2(\tilde x_4)$} [l] at 281 294
\pinlabel {$\tau_2(\tilde y_5)$} [l] at 300 271
\pinlabel {$\tau_2(\tilde y_4)$} [r] at 89 315
\pinlabel {$\tau_2(\tilde x_5)$} [r] at 65 280
\pinlabel $*$ at 117 315
\pinlabel $*$ at 246 315
\pinlabel $*$ at 337 222
\pinlabel $*$ at 337 94
\pinlabel $*$ at 245 2
\pinlabel $*$ at 116 2
\pinlabel $*$ at 25 94
\pinlabel $*$ at 25 222
\endlabellist
\includegraphics[width=3.4in]{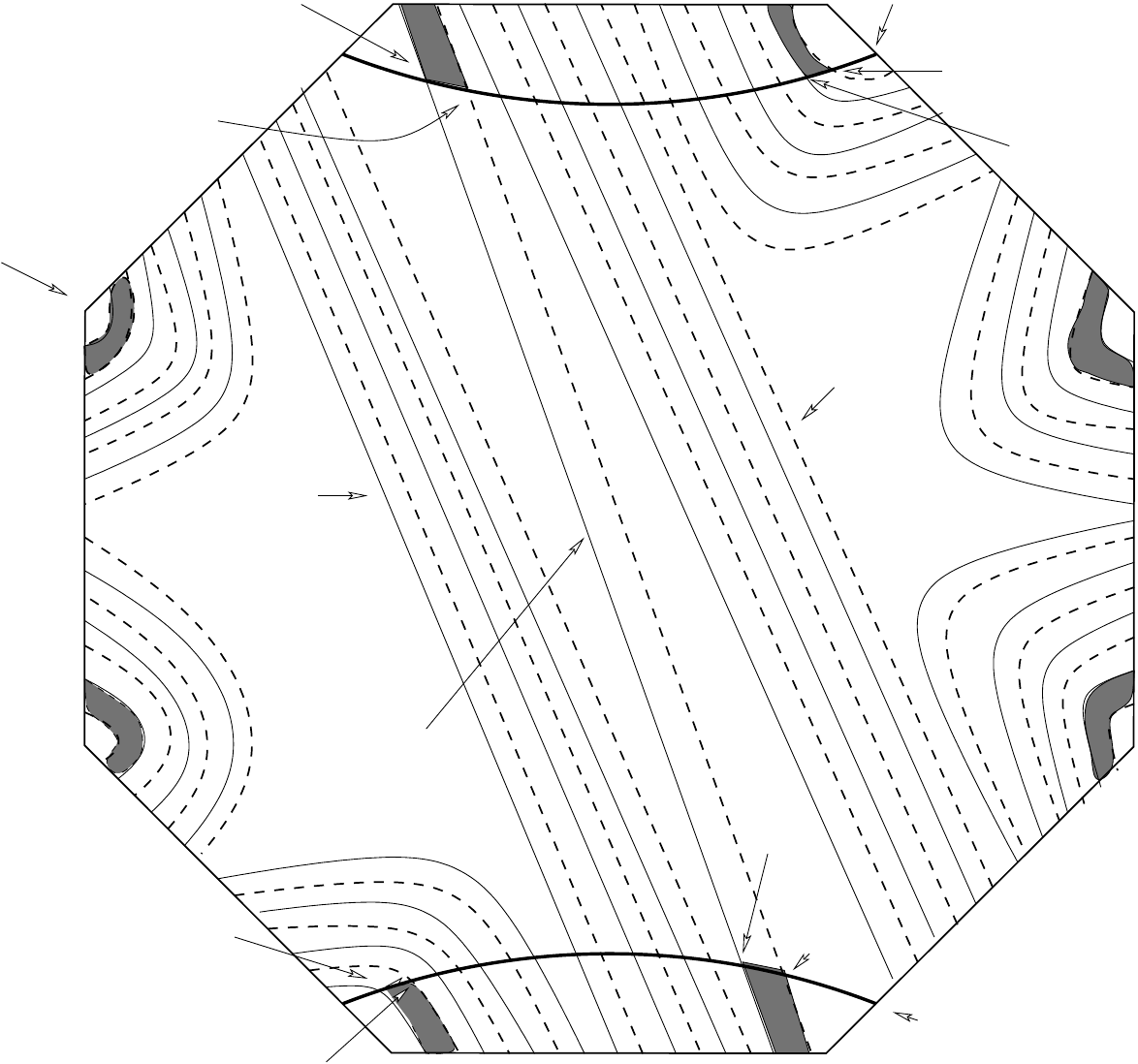}
\end{center}
\caption{Two holomorphic disks, one in $\pi_2((\widetilde{y_5}, \tau_2(\widetilde{y_4})),(\widetilde{x_4}, \tau_2(\widetilde{x_5})))$ and one in $\pi_2((\widetilde{y_4}, \tau_2(\widetilde{y_5})),(\widetilde{x_5}, \tau_2(\widetilde{x_4})))$}
\label{fig:K157Branchedtnon}
\end{figure}

\begin{figure}[ht!]
\begin{center}
\labellist
\pinlabel $x_1$ [b] at 2 116
\pinlabel $x_2$ [b] at 25 116
\pinlabel $x_3$ [b] at 46 116
\pinlabel $x_4$ [b] at 69 116
\pinlabel $x_5$ [b] at 105 116
\pinlabel $x_6$ [b] at 131 116
\pinlabel $x_7$ [b] at 159 116
\pinlabel $x_8$ [b] at 183 116
\pinlabel $y_7$ [t] at 15 95
\pinlabel $y_6$ [t] at 37 95
\pinlabel $y_5$ [t] at 57 95
\pinlabel $y_4$ [t] at 88 95
\pinlabel $y_3$ [t] at 121 95
\pinlabel $y_2$ [t] at 146 95
\pinlabel $y_1$ [t] at 172 95
\pinlabel $w$ [t] at 86 114
\pinlabel $z$ [t] at 117 70
\pinlabel $*$ [b] at 88 123
\pinlabel $*$ [b] at 108 87
\endlabellist
\includegraphics[width=3.3in]{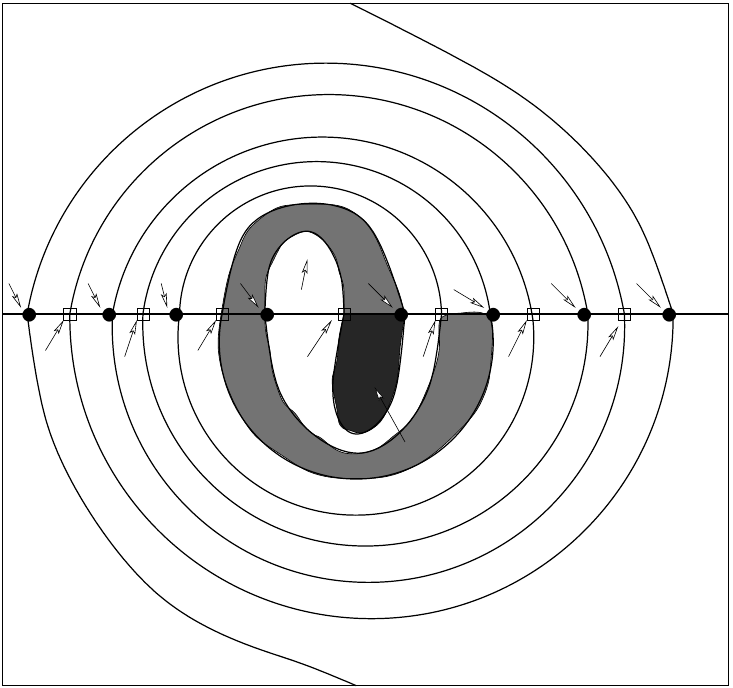}
\end{center}
\caption{Two holomorphic disks, one in $\pi_2(x_5,y_4)$ and the other in $\pi_2(x_6,y_3)$}
\label{fig:K157tnon2}
\end{figure}

\begin{figure}[ht!]
\begin{center}
\labellist
\pinlabel $w$ [r] at 0 237
\pinlabel $*$ [lb] at 174 156
\pinlabel $\beta$ [r] at 93 168
\pinlabel $z$ [rt] at 126 99
\pinlabel $\tau_2(\beta)$ [b] at 247 199
\pinlabel $\tau_2(\alpha)$ [bl] at 265 314
\pinlabel $\alpha$ [l] at 273 10
\pinlabel {$\tilde y_3$} [b] at 242 32
\pinlabel {$\tilde x_6$} [b] at 228 60
\pinlabel {$\tau_2(\tilde y_4)$} [r] at 89 315
\pinlabel {$\tau_2(\tilde x_5)$} [r] at 65 280
\pinlabel $*$ at 118 313
\pinlabel $*$ at 246 313
\pinlabel $*$ at 337 222
\pinlabel $*$ at 337 92
\pinlabel $*$ at 245 1
\pinlabel $*$ at 116 1
\pinlabel $*$ at 25 92
\pinlabel $*$ at 25 222
\endlabellist
\includegraphics[width=3.4in]{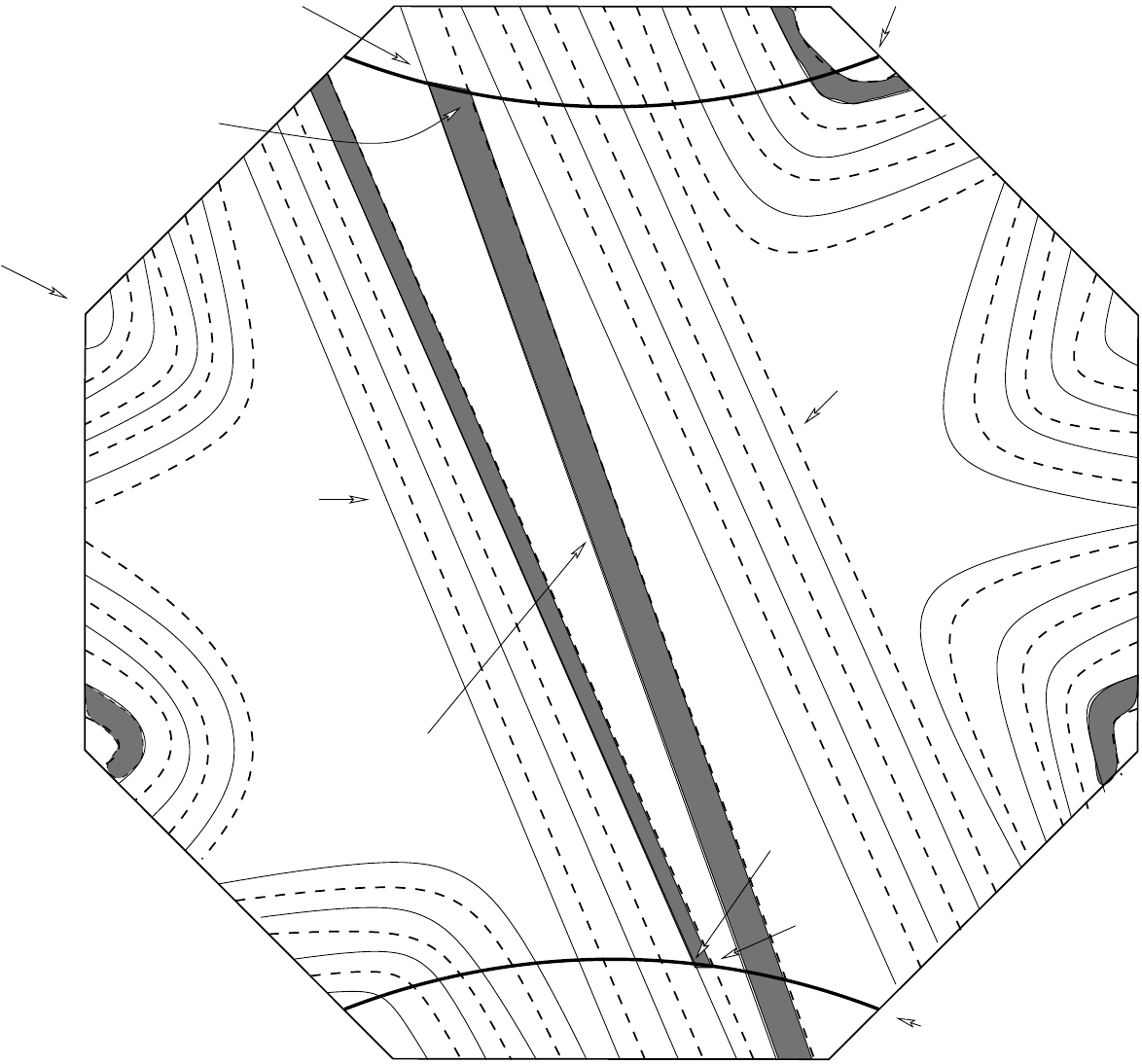}
\end{center}
\caption{A holomorphic disk in $\pi_2((\widetilde{x_6},\tau_2(\widetilde{x_5})), (\tau_2(\widetilde{y_4}), \widetilde{y_3}))$}
\label{fig:K157Branchedtnon2}
\end{figure}

\begin{figure}[ht!]
\begin{center}
\labellist
\pinlabel $(x_7,`x_8)$ [l] at 371 376
\pinlabel $(x_6,`x_7)$ [l] at 371 348
\pinlabel $(x_5,`x_6)$ [l] at 371 319
\pinlabel $(y_4,``y_5)(y_1,``y_2)$ [l] at 149 247
\pinlabel $(y_5,``y_6)(y_2,``y_3)$ [l] at 149 220
\pinlabel $(y_6,``y_7)(y_3,``y_4)$ [l] at 149 193
\pinlabel $(x_1,`x_2)$ [r] at 30 110
\pinlabel $(x_2,`x_3)$ [r] at 30 82
\pinlabel $(x_3,`x_4)$ [r] at 30 54
\pinlabel $(x_4,`x_5)$ [t] at 200 97
\pinlabel {filtration level} [t] at 0 0
\pinlabel ${-}1$ <0pt,-15pt> at  0 0
\pinlabel {filtration level} [t] at 200 0
\pinlabel $0$  <0pt,-15pt> at 200 0
\pinlabel {filtration level} [t] at 400 0
\pinlabel $1$ <0pt,-15pt> at 400 0
\endlabellist
\includegraphics[width=3.1in]{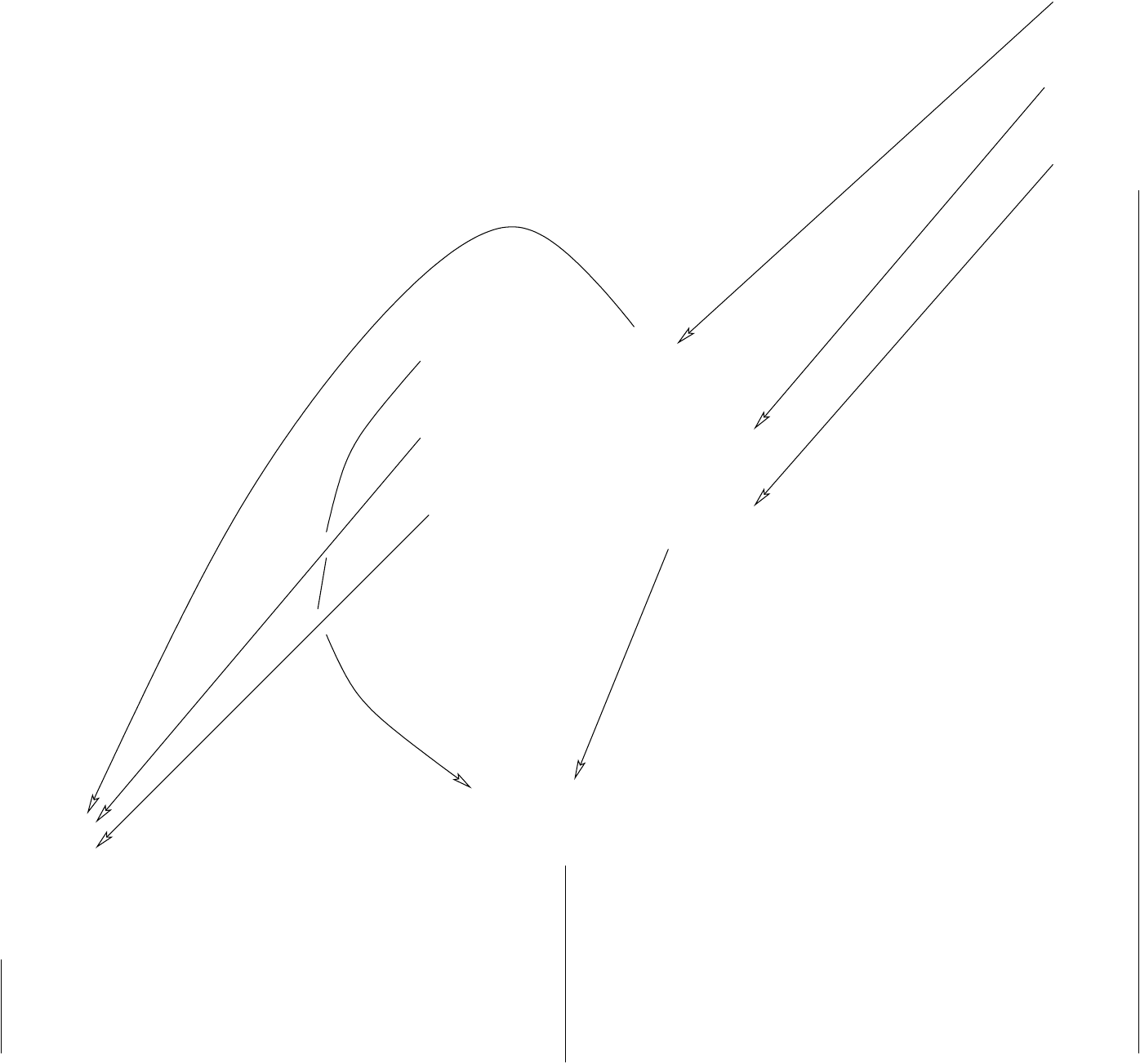}
\end{center}
\vspace{2mm}
\caption{$\widehat{CFK}(\Sigma^2(K), \widetilde{K}(15,7), \mathfrak{s}_{\pm 2})$}
\label{fig:K157Branchedspin2}
\end{figure}

\begin{figure}[ht!]
\begin{center}
\labellist
\pinlabel $x_8$ [t] at -6 68
\pinlabel $x_1$ [t] at 20 68
\pinlabel $x_6$ [t] at 46 68
\pinlabel $x_7$ [t] at 72 68
\pinlabel $x_5$ [b] <2pt,0pt> at 323 265
\pinlabel $x_4$ [b] <2pt,0pt> at 298 265
\pinlabel $x_3$ [b] <2pt,0pt> at 273 265
\pinlabel $x_2$ [b] <2pt,0pt> at 247 265
\pinlabel $y_1$ at 114 170
\pinlabel $y_7$ at 134 170
\pinlabel $y_3$ at 155 170
\pinlabel $y_6$ at 175 170
\pinlabel $y_2$ at 197 170
\pinlabel $y_5$ at 218 170
\pinlabel $y_4$ at 250 170
\pinlabel {filtration level} [t] at 34 0
\pinlabel ${-}1$ <0pt,-15pt> at 34 0
\pinlabel {filtration level} [t] at 177 0
\pinlabel $0$  <0pt,-15pt> at 177 0
\pinlabel {filtration level} [t] <10pt,0pt> at 287 0
\pinlabel $1$ <10pt,-15pt> at 287 0
\endlabellist
\includegraphics[width=2.8in]{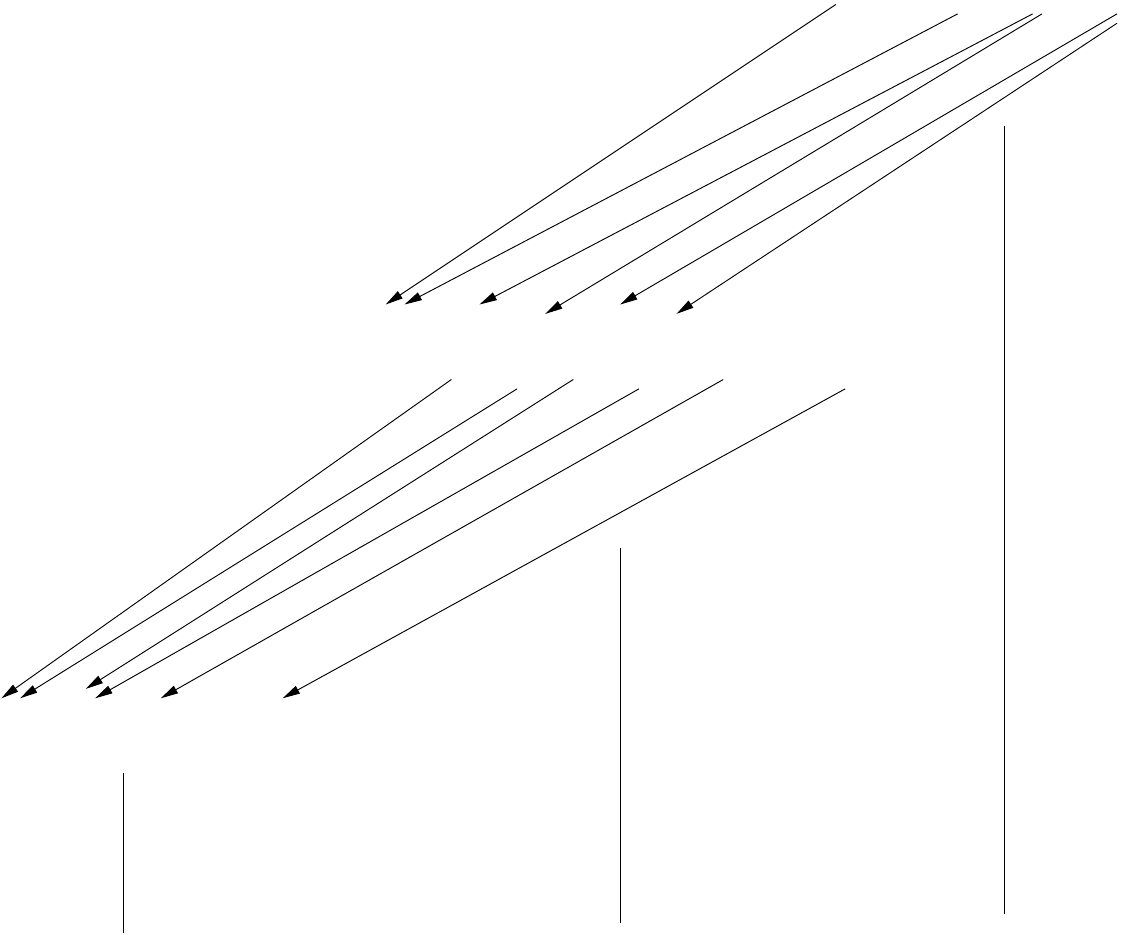}
\end{center}
\vspace{2mm}
\caption{$\mathbb{Z}$--filtered chain complex for $\widehat{CFK}(S^3; K(15,4))$}
\label{fig:K154cfhat}
\end{figure}

$\widehat{HFK}(\Sigma^2(K); \widetilde{K})$ distinguishes between knots with the same knot Floer homology.  In \cite{MR1988285} \Ozsvath and \Szabo prove that for $K$ an alternating knot in $S^3$, the knot Floer homology is determined completely by the Alexander polynomial and the signature.

K(15,7) and K(15,4) are two-bridge knots with the same Alexander polynomial and signature but differing $\mathbb{Z}_2$--graded knot Floer homologies in the double branched cover.\footnote{The double-branched covers of $K(15,7)$ and $K(15,4)$ are different 3-manifolds (the lens spaces -L(15,7) and -L(15,4), respectively).  In fact, $\widehat{HF}(\Sigma^2(K(15,7))) \not\cong \widehat{HF}(\Sigma^2(K(15,4)))$ as $\mathbb{Q}$--graded groups (See Proposition 4.8 in \cite{MR1957829} for an inductive formula for the $\mathbb{Q}$ grading of generators in lens spaces).  This immediately implies that $\widehat{HFK}(\Sigma^2(K(15,7));\widetilde{K}(15,7)) \not\cong \widehat{HFK}(\Sigma^2(K(15,4));\widetilde{K}(15,4))$ as $\mathbb{Q}$--graded groups.  We will show the stronger statement that $\widehat{HFK}(\Sigma^2(K(15,7)); \widetilde{K}(15,7)) \not\cong \widehat{HFK}(\Sigma^2(K(15,4)); \widetilde{K}(15,4))$ as $\mathbb{Z}_2$--graded groups.} The computations of both $\widehat{HFK}(S^3;K)$ and $\widehat{HFK}(\Sigma^2(K);\widetilde{K})$ for these two knots are given below.

\subsubsection*{Computation for $K(15,7)$}$\phantom{99}$

We start by computing $\widehat{HFK}(S^3;K)$ for $K = K(15,7)$.  In the genus $1$ Heegaard diagram compatible with $K$ given in \fullref{fig:K157tinv}, we have $15$ generators, which we label $x_1, \ldots, x_8$ and $y_1, \ldots, y_7$.  This Heegaard diagram for $S^3$ compatible with $K$ was obtained by taking the handlebody decomposition of $S^3 - K$ coming from the Schubert normal form for $K$ and destabilizing once.

The differential counts maps of the disk into $\Sym^1(S) = S$.  See \fullref{fig:K157tinv} for an example.

The filtered chain complex for $\widehat{CF}$ is pictured in \fullref{fig:K157cfhat}.  The $i$-th vertical slice of this filtered chain complex is the chain complex $\mathcal{F}_i/\mathcal{F}_{i-1} = \widehat{CFK}(S^3;K,i)$.

We construct a Heegaard diagram for $\Sigma^2(K)$ compatible with $\widetilde{K}$ by taking the branched double cover of the Heegaard surface $\Sigma$ around the two basepoints $w$ and $z$.  The $\alpha$ and $\beta$ curves of our original Heegaard diagram then lift to two $\alpha$ curves and two $\beta$ curves, and $w$ and $z$ lift to the two basepoints for the doubly-pointed Heegaard diagram for $\Sigma^2(K)$.

All generators of $\widehat{HFK}(\Sigma^2(K); \widetilde{K})$ are of the form $(\widetilde{x_i}, \tau_2(\widetilde{x_j}))$ or $(\widetilde{y_i}, \tau_2(\widetilde{y_j}))$, where $\tau_2$ is the non-trivial element of $\mathbb{Z}_2$.  Under the $\mathbb{Z}_2$ action on the $\alpha$ and $\beta$ curves, we get a natural $\mathbb{Z}_2$ action on these generators: $$\tau_2(\widetilde{x_i}, \tau_2(\widetilde{x_j})) = (\tau_2(\widetilde{x_i}), \widetilde{x_j}).$$  We shall refer to such a pair of generators as a \begin{it} conjugate pair \end{it}, since they are in conjugate Spin$^c$ structures with respect to the central Spin$^c$ structure, $\mathfrak{s}_0$.

\begin{figure}[ht!]
\begin{center}
\labellist
\pinlabel $x_1$ [b] at 47 291
\pinlabel $x_7$ [b] at 118 281
\pinlabel $x_5$ [b] at 200 287
\pinlabel $x_3$ [b] at 263 287
\pinlabel $x_8$ [b] at 322 291
\pinlabel $x_6$ [b] at 393 299
\pinlabel $x_4$ [b] at 468 295
\pinlabel $x_2$ [b] at 523 281
\pinlabel $y_2$ [t] at 80 241
\pinlabel $y_4$ [t] at 158 241
\pinlabel $y_6$ [t] at 235 241
\pinlabel $y_1$ [t] at 303 241
\pinlabel $y_3$ [t] at 366 241
\pinlabel $y_5$ [t] at 452 241
\pinlabel $y_7$ [t] at 501 227
\pinlabel $w$ [b] at 155 313
\pinlabel $*$ [t] at 159 296
\pinlabel* $z$ [br] at 392 252
\pinlabel* $*$ [tl] at 402 241
\endlabellist
\includegraphics[width=3.5in]{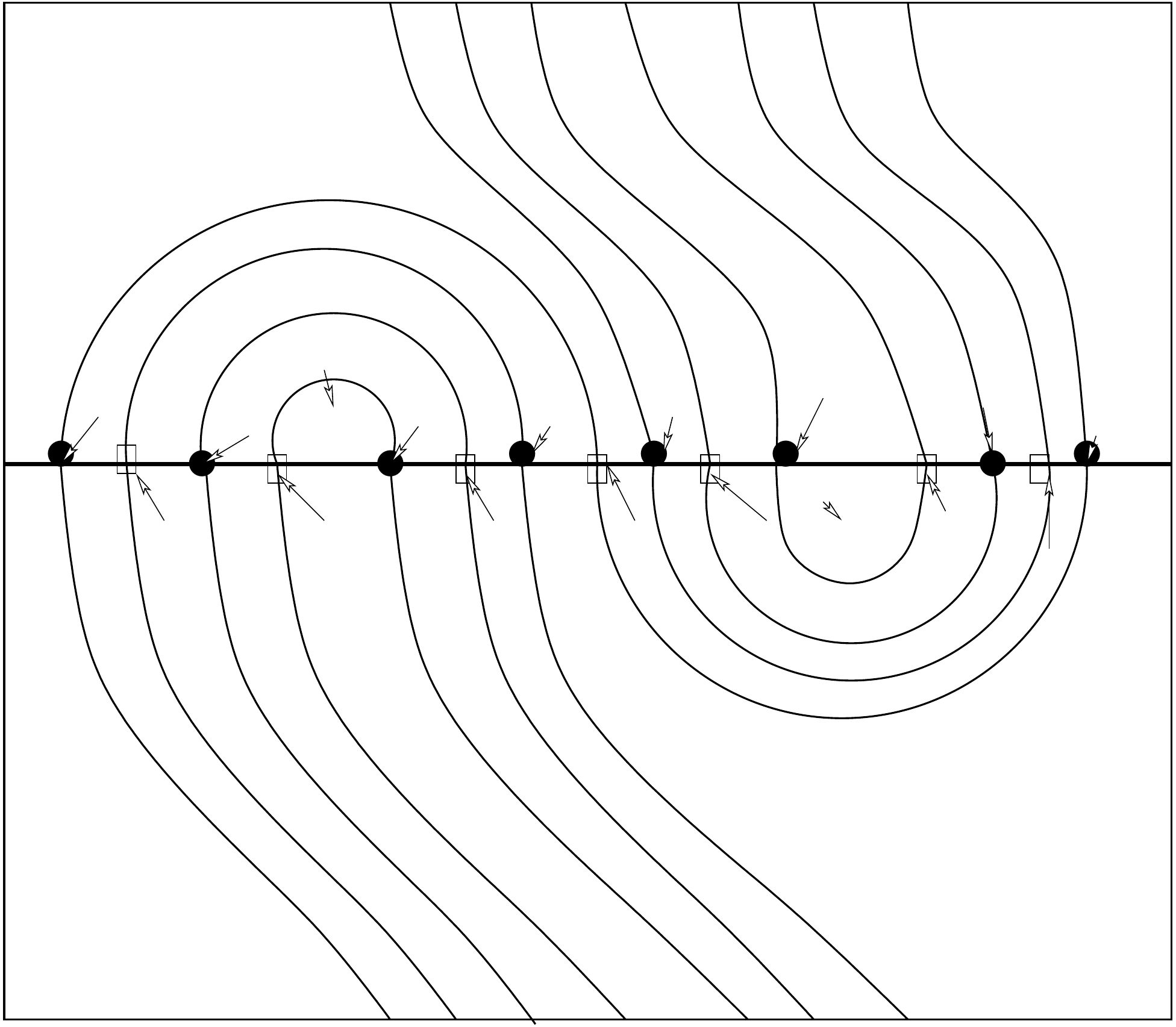}
\end{center}
\caption{Genus $1$ (destabilized) Heegaard diagram for $K(15,4)$ in $S^3$}
\label{fig:K154}
\end{figure}

\begin{figure}[ht!]
\begin{center}
\labellist
\pinlabel $w$ [r] at 0 391
\pinlabel $\beta$ [r] at 223 186
\pinlabel $z$ [lb] at 398 351
\pinlabel $*$ [tr] at 330 303
\pinlabel $\tau_2(\beta)$ [bl] at 405 234
\pinlabel $\tau_2(\alpha)$ [l] at 525 548
\pinlabel $\alpha$ [tl] at 509 0
\pinlabel $*$ at 214 567
\pinlabel $*$ at 450 567
\pinlabel $*$ at 615 400
\pinlabel $*$ at 615 167
\pinlabel $*$ at 448 0
\pinlabel $*$ at 214 0
\pinlabel $*$ at 48 166
\pinlabel $*$ at 48 402
\endlabellist
\includegraphics[width=3.4in]{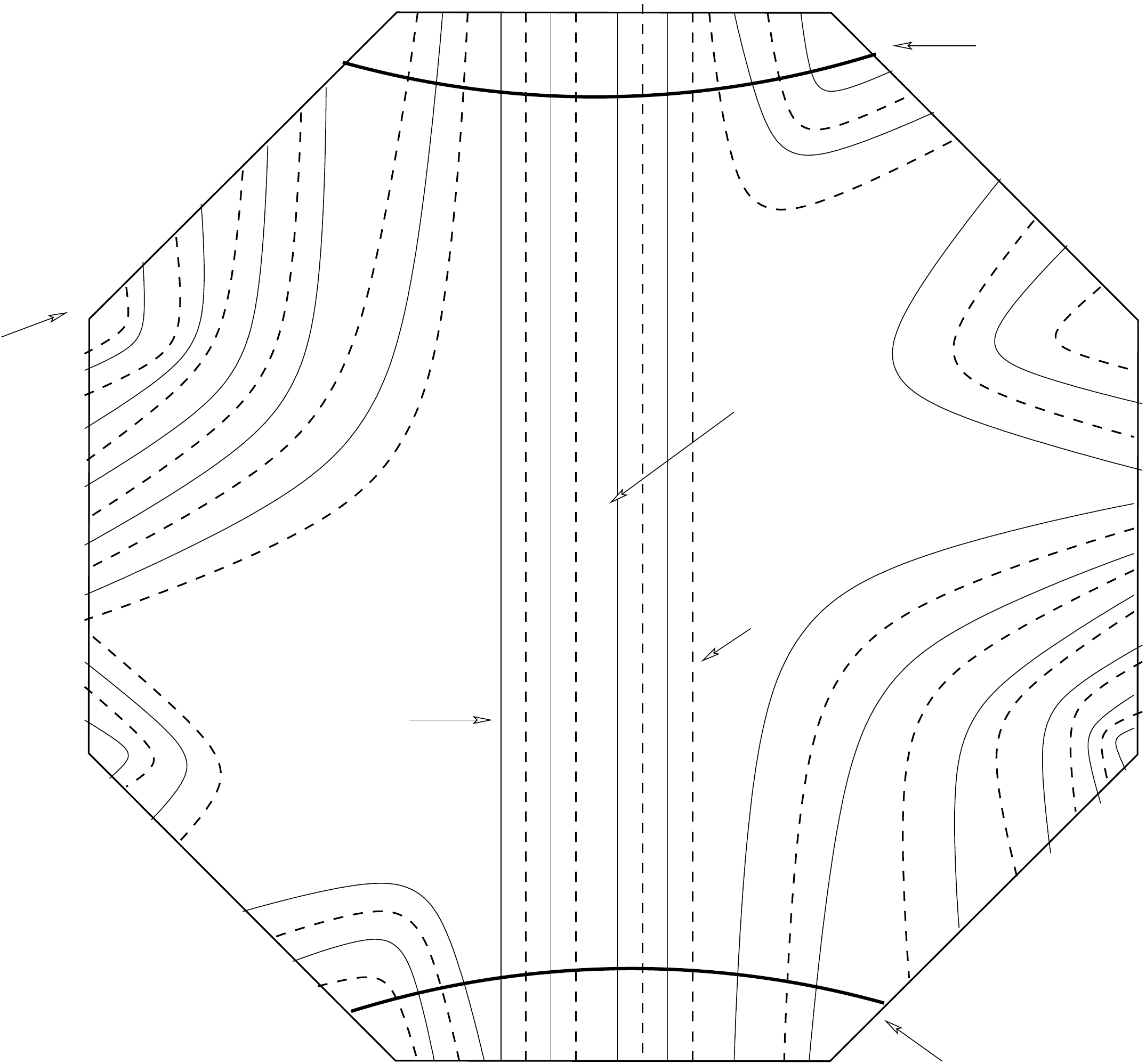}
\end{center}
\caption{Genus $2$ (destabilized) Heegaard diagram for $\widetilde{K}(15,4)$ in $L(15,4)$}
\label{fig:K154Branched}
\end{figure}

{\small
\begin{table}[ht!]
\centering
\begin{tabular}{|c|c|c|c|c|}
\hline
$\mathfrak{s}_0$ & $\mathfrak{s}_{\pm 1}$ & $\mathfrak{s}_{\pm 2}$ & $\mathfrak{s}_{\pm 3}$ & $\mathfrak{s}_{\pm 4}$ \\[-0.5pt]
\hline
$(x_1,x_1), (y_1,y_1)$ & $(x_2,x_3), (y_2,y_3)$ &  $(x_1,x_2)$ & $(x_3,x_4), (y_3,y_4)$ & $(x_2,x_4), (y_2,y_4)$ \\[-0.5pt]
$(x_2,x_2), (y_2,y_2)$ & $(x_4,x_5), (y_4,y_5)$ &              & $(x_5,x_6), (y_5,y_6)$ & $(x_3,x_5), (y_3,y_5)$ \\[-0.5pt]
$(x_3,x_3), (y_3,y_3)$ & $(x_6,x_7), (y_6,y_7)$ &              & $(x_7,x_8), (y_1,y_2)$ & $(x_4,x_6), (y_4,y_6)$ \\[-0.5pt]
$(x_4,x_4), (y_4,y_4)$ & $(x_1,x_8)$            &              & $(x_2,x_8), (y_1,y_7)$ & $(x_5,x_7), (y_5,y_7)$ \\[-0.5pt]
$(x_5,x_5), (y_5,y_5)$ &                        &              & $(x_1,x_3)$            & $(x_6,x_8), (y_1,y_3)$ \\[-0.5pt]
$(x_6,x_6), (y_6,y_6)$ &                        &              &                        & $(x_3,x_8), (y_1,y_6)$ \\[-0.5pt]
$(x_7,x_7), (y_7,y_7)$ &                        &              &                        & $(x_1,x_7)$            \\[-0.5pt]
$(x_8,x_8)$            &                        &              &                        &                        \\[-0.5pt]
\hline
\end{tabular}

\vspace{1mm}
\begin{tabular}{|c|c|c|}
\hline
$\mathfrak{s}_{\pm 5}$ & $\mathfrak{s}_{\pm 6}$ & $\mathfrak{s}_{\pm 7}$ \\
\hline
$(x_2,x_5),(y_2,y_5)$ & $(x_2,x_7),(y_2,y_7)$ & $(x_1,x_5), (y_1,y_5)$\\
$(x_4,x_7),(y_4,y_7)$ & $(x_1,x_4)$           & $(x_2,x_6), (y_2,y_6)$\\
$(x_1,x_6)$           &                       & $(x_3,x_7), (y_3,y_7)$\\
                      &                       & $(x_3,x_6), (y_3,y_6)$\\
                      &                       & $(x_5,x_8), (y_1,y_4)$\\
                      &                       & $(x_4,x_8)$\\
\hline
\end{tabular}
\vspace{1mm}
\caption{Spin$^c$ structures $\mathfrak{s}_{0}, \ldots, \mathfrak{s}_{\pm 7}$ for $\Sigma^2(K(15,4))$}
\label{table:Spinc154}
\end{table}
}

\begin{figure}[ht!]
\vspace{4mm}
\begin{center}
\labellist
\pinlabel $(x_2,`x_3)$ [l] <0pt, 12pt> at 163 137
\pinlabel $(x_4,`x_5)$ [l] at 163 137
\pinlabel $(x_6,`x_7)$ [r] <0pt, 12pt> at 16 16
\pinlabel $(x_1,`x_8)$ [r] at 16 16
\pinlabel $(y_4,``y_5)$ [l] at 73 84
\pinlabel $(y_2,``y_3)$ [l] <0pt, -12pt> at 73 84
\pinlabel $(y_6,``y_7)$ [l] <0pt, 12pt> at 73 84
\pinlabel {filtration level} [t] at 0 0
\pinlabel ${-}1$ <0pt,-15pt> at  0 0
\pinlabel {filtration level} [t] at 86 0
\pinlabel $0$  <0pt,-15pt> at 86 0
\pinlabel {filtration level} [t] at 177 0
\pinlabel $1$ <0pt,-15pt> at 177 0
\endlabellist
\includegraphics[width=3in]{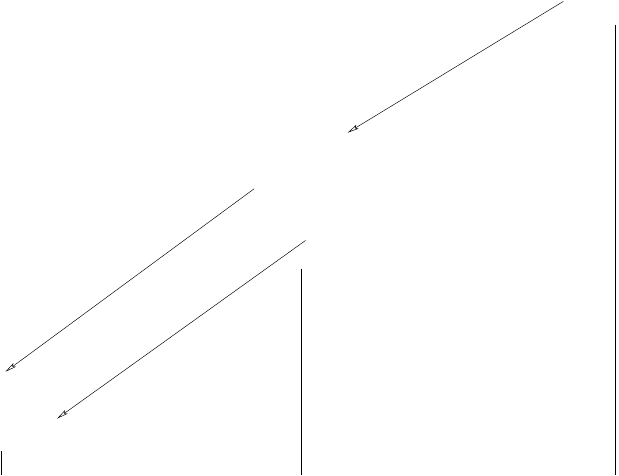}
\end{center}
\vspace{2mm}
\caption{$\widehat{CFK}(\Sigma^2(K), \widetilde{K}(15,4), \mathfrak{s}_{\pm 1})$}
\label{fig:K154Branchedspin1}
\end{figure}

\begin{figure}[ht!]
\vspace{4mm}
\begin{center}
\labellist
\pinlabel $(x_1,`x_2)$ [b] at 0 26
\pinlabel {filtration level} [t] at 0 0
\pinlabel $0$  <0pt,-15pt> at 0 0
\endlabellist
\includegraphics[height=0.5in]{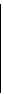}
\end{center}
\vspace{2mm}
\caption{$\widehat{CFK}(\Sigma^2(K), \widetilde{K}(15,4), \mathfrak{s}_{\pm 2})$}
\label{fig:K154Branchedspin2}
\end{figure}

\begin{figure}[ht!]
\begin{center}
\labellist
\pinlabel $(x_2,`x_4)$ [l] <0pt, -12pt> at 177 158
\pinlabel $(x_3,`x_5)$ [l] at 177 158
\pinlabel $(x_5,`x_7)$ [l] <0pt, -24pt> at 87 167
\pinlabel $(x_4,`x_6)$ [l] <0pt, -12pt> at 87 167
\pinlabel $(x_3,`x_8)$ [l] at 87 167
\pinlabel $(x_6,`x_8)$ [r] at 15 24
\pinlabel $(x_1,`x_7)$ [r] <0pt, -12pt> at 15 24

\pinlabel $(y_1,``y_3)$ [l] at 56 76
\pinlabel $(y_1,``y_6)$ [l] <0pt, 12pt> at 56 76
\pinlabel $(y_5,``y_7)$ [l] at 89 76
\pinlabel $(y_3,``y_5)$ [l] <0pt, 12pt> at 89 76
\pinlabel $(y_2,``y_4)$ [r] at 142 76
\pinlabel $(y_4,``y_6)$ [r] <0pt, 12pt> at 142 76

\pinlabel {filtration level} [t] at 0 0
\pinlabel ${-}1$ <0pt,-15pt> at  0 0
\pinlabel {filtration level} [t] at 100 0
\pinlabel $0$  <0pt,-15pt> at 100 0
\pinlabel {filtration level} [t] at 193 0
\pinlabel $1$ <0pt,-15pt> at 193 0
\endlabellist
\includegraphics[width=3.5in]{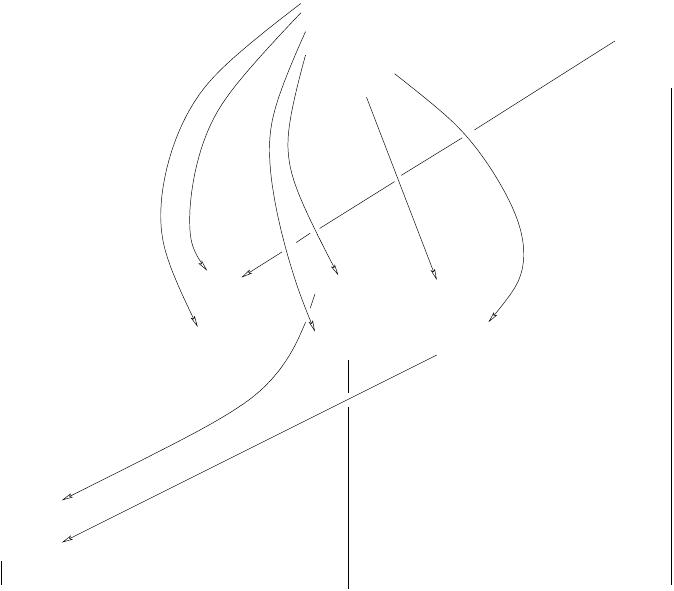}
\end{center}
\vspace{2mm}
\caption{$\widehat{CFK}(\Sigma^2(K), \widetilde{K}(15,4), \mathfrak{s}_{\pm 4})$}
\label{fig:K154Branchedspin4}
\end{figure}

\begin{figure}[ht!]
\vspace{4mm}
\begin{center}
\labellist
\pinlabel $(y_3,``y_7)$ [b] at 36 132
\pinlabel $(y_3,``y_6)$ [b] at 76 132
\pinlabel $(y_2,``y_6)$ [b] at 116 132
\pinlabel $(x_2,`x_6)$ [b] at 0 80
\pinlabel $(x_4,`x_8)$ [b] at 30 80
\pinlabel $(x_3,`x_6)$ [b] at 60 80
\pinlabel $(x_5,`x_8)$ [b] at 90 80
\pinlabel $(x_3,`x_7)$ [b] at 120 80
\pinlabel $(x_1,`x_5)$ [b] at 150 80
\pinlabel $(y_1,``y_5)$ [t] at 53 38
\pinlabel $(y_1,``y_4)$ [t] at 101 38
\pinlabel {filtration level} [t] at 77 0
\pinlabel $0$  <0pt,-15pt> at 77 0
\endlabellist
\includegraphics[width=2.5in]{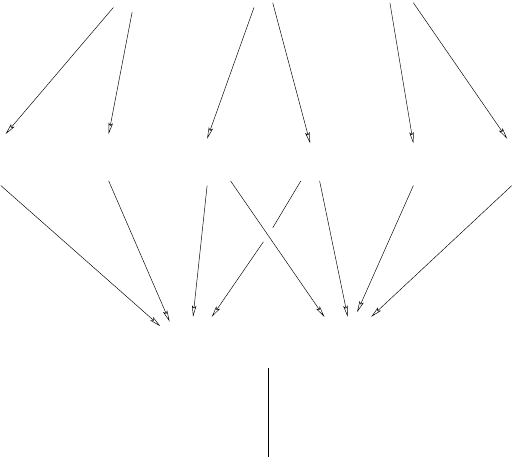}
\end{center}
\vspace{2mm}
\caption{$\widehat{CFK}(\Sigma^2(K), \widetilde{K}(15,4), \mathfrak{s}_{\pm 7})$}
\label{fig:K154Branchedspin7}
\end{figure}

\subsubsection*{$\tau_2$--invariant generators}$\phantom{99}$

The $\tau_2$--invariant generators are those of the form $(\widetilde{x_i}, \tau_2(\widetilde{x_i}))$ or $(\widetilde{y_i}, \tau_2(\widetilde{y_i}))$.  Furthermore, for every topological bigon downstairs we see a corresponding topological quadrilateral upstairs (as detailed in the proof of \fullref{lemma:mgrading}).  See the shaded disk in $\pi_2(y_5,x_4)$ in \fullref{fig:K157tinv} and its lift in $\pi_2((\widetilde{y}_5, \tau_2(\widetilde{y}_5)), (\widetilde{x}_4, \tau_2(\widetilde{x}_4)))$ in \fullref{fig:K157Branchedtinv}.

\subsubsection*{$\tau_2$--non-invariant generators}$\phantom{99}$

We will first state how the generators split up according to Spin$^c$ structures, then explicity compute $\widehat{HFK}(\Sigma^2(K); \widetilde{K}; \mathfrak{s}_{\pm 2})$.  Here we use $\mathfrak{s}_{\pm 2}$ to denote the conjugate Spin$^c$ structures on $\Sigma^2(K)$ corresponding to $\pm 2 \in H_1(\Sigma^2(K)) \cong \mathbb{Z}_{15}$.

Each column in Table \ref{table:Spinc157} gives the generators in two conjugate Spin$^c$ structures simultaneously (except in the case of the ``central'' Spin$^c$ structure, $\mathfrak{s}_0$, which is its own conjugate).  When we write $(x_1, x_8)$, for example, we refer to two generators simultaneously: $(\widetilde{x}_1, \tau_2(\widetilde{x}_8)) \in \mathfrak{s}_1$ and $(\tau_2(\widetilde{x}_1), \widetilde{x}_8) \in \mathfrak{s}_{-1}$.

Now let's more closely examine one of the Spin$^c$ structures, $\mathfrak{s}_{\pm 2}$.

To see that these are the generators in the two conjugate Spin$^c$ structures $\mathfrak{s}_{\pm 2}$, connect, for example, $(\widetilde{x}_1, \tau_2(\widetilde{x}_1))$ to $(\widetilde{x}_1, \tau_2(\widetilde{x}_2))$ by the path in $\Sym^2(\widetilde{S})$ which is the product of the constant path $\widetilde{x}_1 \rightarrow \widetilde{x}_1$ with the path $\tau_2(\widetilde{x}_1) \rightarrow \tau_2(\widetilde{x_2})$ along the $\tau_2(\widetilde{\beta})$ curve.  Close the path to a loop $\gamma$ by taking a path $\tau_2(\widetilde{x}_1) \rightarrow \tau_2(\widetilde{x}_2)$ along the $\tau_2(\widetilde{\alpha})$ curve.  Since $$\widetilde{\alpha} \cdot \widetilde{\beta} = \tau_2(\widetilde{\alpha}) \cdot \tau_2(\widetilde{\beta}) = 1,$$ and $$\widetilde{\alpha} \cdot \tau_2(\widetilde{\beta}) = \tau_2(\widetilde{\alpha}) \cdot \widetilde{\beta} = -1,$$ we see that if we assert that $\widetilde{\alpha}$ is the positive generator of $H_1(\Sigma^2(K))$, then the loop $\gamma$ represents the element $2 \in H_1(\Sigma^2(K))$.  

We see this because a pushoff of $\gamma$ has one intersection with $\tau_2(\widetilde{\alpha})$ with local multiplicity $-1$ (at $\tau_2(\widetilde{y_1})$) and one intersection with $\widetilde{\alpha}$ with local multiplicity $1$ (at $\widetilde{x_2}$).  Recall that $\tau_2(\widetilde{\alpha}) = -\widetilde{\alpha}$ in $H_1(\Sigma^2(K))$.

The same type of calculation can be performed to verify that all of the other generators are in the stated Spin$^c$ structures.

Now focus on, for example, $(\widetilde{x_4}, \tau_2(\widetilde{x_5}))$ and $(\widetilde{y_5}, \tau_2(\widetilde{y_4}))$ (and the corresponding pair $(\tau_2(\widetilde{x_4}), \widetilde{x_5}), (\tau_2(\widetilde{y_5}), \widetilde{y_4})$ in the conjugate Spin$^c$ structure).

Notice that in the Heegaard diagram for $(S^3,K(15,7))$ pictured in \fullref{fig:K157tnon} we see a topological disk in $\pi_2(x_5, y_4)$ and a topological disk in $\pi_2(y_5, x_4)$, both with $n_w = 0$ and $n_z = 1$.  

If we lift both of these to the Heegaard diagram for $(\Sigma^2(K); \widetilde{K})$ we get two quadrilaterals, one representing a disk in $\pi_2((\widetilde{x_5}, \tau_2(\widetilde{x_5})),(\widetilde{y_4}, \tau_2(\widetilde{y_4})))$ and one in $$\pi_2((\widetilde{y_5}, \tau_2(\widetilde{y_5})), (\widetilde{x_4}, \tau_2(\widetilde{x_4}))).$$
Taking the difference of these two disks yields two disks, one in $$\pi_2((\widetilde{y_4}, \tau_2(\widetilde{y_5})),(\widetilde{x_5}, \tau_2(\widetilde{x_4})))$$ and one in $\pi_2((\widetilde{y_5}, \tau_2(\widetilde{y_4})),(\widetilde{x_4}, \tau_2(\widetilde{x_5}))).$  Both of these disks have $n_z = n_w = 0$ and Maslov index $1$.  Both of these disks are holomorphic because they are topological quadrilaterals (see \fullref{fig:K157Branchedtnon}).

A similar argument shows that there is another pair of disks, one in $$\pi_2((\widetilde{y}_4, \tau_2(\widetilde{y}_3)),(\widetilde{x}_4, \tau_2(\widetilde{x}_5)))$$ and one in $$\pi_2((\tau_2(\widetilde{y}_4), \widetilde{y}_3),(\tau_2(\widetilde{x}_4), \widetilde{x}_5)),$$ both with $n_z = n_w = 0$.

There also exist topological disks in $\Sym^2(\widetilde{S})$ which are the sum of a pair of lifts of disks.  Focus, for example, on the pair $(\widetilde{x}_6, \tau_2(\widetilde{x}_5))$ and $(\tau_2(\widetilde{y}_4),\widetilde{y}_3)$.  In the Heegaard diagram downstairs, we see a disk in $\pi_2(x_6, y_3)$ with $n_w = 0$ and $n_z = 1$.  We also see a disk in $\pi_2(x_5, y_4)$.  See \fullref{fig:K157tnon2}.

The sum of the lifts of these two disks to the Heegaard diagram for $\Sigma^2(K)$ again breaks up into two disks.  One of the disks is in $$\pi_2((\widetilde{x_6},\tau_2(\widetilde{x_5})), (\tau_2(\widetilde{y_4}), \widetilde{y_3}))$$ and the other is between the conjugate generators, ie, in $$\pi_2((\tau_2(\widetilde{x_6}), \widetilde{x_5}), (\widetilde{y_4}, \tau_2(\widetilde{y_3}))).$$  Both of these disks have $n_w = 0$ and $n_z = 1$.  See \fullref{fig:K157Branchedtnon2} for one of the two disks.

Using similarly obtained disks, we easily calculate the relative filtration levels and Maslov grading of all generators in $\mathfrak{s}_{\pm 2}$.  The generators of the $\widehat{CFK}$ complex for $K(15,7)$ as well as arrows corresponding to some of the differentials are pictured in \fullref{fig:K157Branchedspin2}.  We omit the $\sim$'s and the $\tau_2$'s, since we are thinking of this as the chain complex for both $\mathfrak{s}_2$ and $\mathfrak{s}_{-2}$.

Note that we have made no claims about whether there are any more holomorphic disks than the ones described.\footnote{In fact, we cannot possibly have listed all of the $d_1$ differentials, for $(d_0 + d_1)^2 \neq 0$ in \fullref{fig:K157Branchedspin2}.}  However, just based on the information contained in \fullref{fig:K157Branchedspin2} we can see that $\widehat{HFK}(\Sigma^2(K); \widetilde{K}; \mathfrak{s}_{\pm 2})$ must have support in three different filtration levels, with $$\widehat{HFK}(\Sigma^2(K); \widetilde{K}, \mathfrak{s}_{\pm 2}, i-1) \cong \widehat{HFK}(\Sigma^2(K); \widetilde{K}, \mathfrak{s}_{\pm 2}, i+1) \cong \mathbb{Z}_2^3.$$

\eject

\subsubsection*{Computation for $K(15,4)$}$\phantom{99}$

We now turn to computations of $\widehat{HFK}(S^3;K)$ and $\widehat{HFK}(\Sigma^2(K); \widetilde{K})$ for $K = K(15,4)$.  Recall that our aim is to prove that $$\widehat{HFK}(\Sigma^2(K); \widetilde{K}(15,7)) \not\cong \widehat{HFK}(\Sigma^2(K); \widetilde{K}(15,4))$$ as $\mathbb{Z}_2$--graded groups.  We will do so by showing that $\widehat{HFK}(\Sigma^2(K); \widetilde{K}(15,4), \mathfrak{s}_k)$ does not look like $\widehat{HFK}(\Sigma^2(K); \widetilde{K}(15,7), \mathfrak{s}_{\pm 2})$ (computed in the previous subsection) for any $k \in \mathbb{Z}_{15}$ with order $15$.

In the genus $1$ Heegaard diagram compatible with $K$ given in \fullref{fig:K154}, we have $15$ generators, which we label $x_1, \ldots, x_8$ and $y_1, \ldots, y_7$.  See \fullref{fig:K154cfhat} for the filtered chain complex.

We construct a genus $2$ Heegaard diagram for $\Sigma^2(K)$ compatible with $\widetilde{K}$ in exactly the same way we did before (see \fullref{fig:K154Branched}).

\subsubsection*{$\tau_2$--non-invariant generators}$\phantom{99}$

The generators split up according to Spin$^c$ structures as detailed in Table \ref{table:Spinc154}.

Using the methods described in the previous subsection, we find relative filtration and Maslov gradings for $\mathfrak{s}_{\pm 1}, \mathfrak{s}_{\pm 2}, \mathfrak{s}_{\pm 4}$, and $\mathfrak{s}_{\pm 7}$.  The generators of the $\widehat{CFK}$ complex for $K(15,4)$ in these Spin$^c$ structures as well as arrows corresponding to some of the differentials are pictured in Figures \ref{fig:K154Branchedspin1} - \ref{fig:K154Branchedspin7}.

Again, we make no claims about whether there are any more holomorphic disks than the ones described.  However, we have enough information about the homology of the associated graded chain complexes to determine that $$\widehat{HFK}(\Sigma^2(K); \widetilde{K}(15,7), \mathfrak{s}_{\pm 2}) \not\cong \widehat{HFK}(\Sigma^2(K); \widetilde{K}(15,4), \mathfrak{s}_k)$$ 
for any $k$ relatively prime to $15$.

\bibliographystyle{gtart}
\bibliography{link}

\end{document}